\documentclass[12pt]{amsart}
\usepackage{lipsum}
\usepackage{amssymb}
\usepackage{enumerate}
\usepackage{amsthm}
\usepackage{tikz-cd}
\usepackage{amsmath}
\usepackage[mathscr]{euscript}
\usepackage[utf8]{inputenc}
\usepackage[english]{babel}
\usepackage{hyperref}
\hypersetup{
    colorlinks=true,
    linkcolor=blue,
    citecolor=red,
    filecolor=red,      
    urlcolor=violet,
}

\makeatletter
\@namedef{subjclassname@2020}{%
  \textup{2020} Mathematics Subject Classification}
\makeatother

\newtheorem{thm}{Theorem}[section]
\newtheorem{lem}[thm]{Lemma}
\newtheorem{prop}[thm]{Proposition}

\newtheorem{corl}[thm]{Corollary}
\newtheorem{xrem}{Remark}
\newtheorem{exm}[thm]{Example}

\DeclareMathOperator{\Nef}{{Nef}}

\DeclareMathOperator{\rk}{{rk}}
\DeclareMathOperator{\Mov}{{Mov}}

\DeclareMathOperator{\Pic}{{Pic}}

\DeclareMathOperator{\Eff}{{Eff}}
\DeclareMathOperator{\End}{{End}}
\DeclareMathOperator{\Div}{{Div}}

\DeclareMathOperator{\Bl}{{Bl}}
\DeclareMathOperator{\NE}{{NE}}

\DeclareMathOperator{\Gr}{{Gr}}

\begin{document}
\baselineskip=18pt

\author{Snehajit Misra}
\author{Nabanita Ray}

\address{Chennai Mathematical Institute, H1 SIPCOT IT Park, Siruseri, Kelambakkam 603103, India.}
\email[Snehajit Misra]{snehajitm@cmi.ac.in, misra08@gmail.com}

\address{Chennai Mathematical Institute, H1 SIPCOT IT Park, Siruseri, Kelambakkam 603103, India.}
\email[Nabanita Ray]{nabanitar@cmi.ac.in}

\begin{abstract}
  Let $E$ be a vector bundle of rank $r$ on a smooth complex projective variety $X$. In this article, we compute the nef and pseudoeffective cones of divisors in the Grassmann bundle $\Gr_X(k,E)$ parametrizing $k$-dimensional subspaces of the fibers of $E$, where $1\leq k \leq \rk(E)$, under assumptions on $X$ as well as on the vector bundle $E$. In particular, we show that nef cone and the pseudoeffective cone of $\Gr_X(k,E)$ coincide if and only if $E$ is a slope semistable bundle on $X$ with $c_2(\End(E))=0$.
  We also discuss about the nefness and ampleness of the universal quotient bundle $Q_k$ on $\Gr_X(k,E)$.
\end{abstract}

\title{Slope Semistability and Positive cones of Grassmann bundles} 
\maketitle
\section{Introduction}
In the last few decades, a number of notions of positivity have been introduced to understand the geometry of the higher dimensional projective varieties. 
 The cone of nef divisors, denoted by $\Nef^1(X)$ and the closure of the cone generated by effective divisors, denoted by $\overline{\Eff}^1(X)$ are two  fundamental invariants of a smooth irreducible projective variety $X$  
  which play a very crucial role in this understanding. The interior of $\Nef^1(X)$ is the ample cone of $X$ which gives important information about the embeddings of $X$ in projective spaces. Also, knowledge of
these cones can be used to study positivity questions, interpolations problems, Seshadri constants. Most of the developments in this direction has been successfully summarized in  \cite{L1}, \cite{L2}.

 In his paper \cite{Miy87}, Yoichi Miyaoka initiated the study of nef and
effective divisors on a projective bundle $\mathbb{P}_C(E)$ over a smooth irreducible curve $C$ in characteristics 0, where $E$ is any rank $r$ vector bundle on $C$, and gives a numerical criterion for slope
semi-stability of $E$ in terms of nefness of the normalized hyperplane class $\lambda_E$ on $\mathbb{P}_C(E)$. More generally, in \cite{Fu}, the cone of effective $k$-cycles in $\mathbb{P}_C(E)$ 
over an irreducible curve $C$ is described in terms of the numerical 
data appearing in  Harder-Narasimhan filtration of the bundle $E$, which generalizes Miyaoka's result in the semistable case. Later  the  nef and pseudo-effective cones of the Grassmannian bundle $\Gr_C(k,E)$ 
parametrizing all the $k$-dimensional quotients of the fibres of the rank $r$ vector bundle $E$ over a smooth curve $C$ defined over an algebraically closed field  of positive characteristic, where $1 \leq k \leq r-1$, have been studied in \cite{B-P} and \cite{B-H-P} respectively. The results in \cite{B-P} and \cite{B-H-P} generalize the known results about nef cone and pseudoeffective cones of projective bundles. However, in most of these cases, the Picard number of these parametrizing spaces is 2, and hence the pseudo-effective cones and nef cones of divisors are finite polyhedra (see \cite{Fu}) generated by two extremal rays in
a two dimensional space. When the Picard number is at least 3, there are very few
examples where  the nef and pseudo-effective cones of divisors are computed. 

In this article, our first main result is the following.
\begin{thm}\label{thm1.1}
 Let $E$ be a semistable vector bundle of rank $r$ on a smooth irreducible complex projective variety $X$ such that $c_2\bigl(\End(E)\bigr) = 2rc_2(E)-(r-1)c^2_1(E) = 0 \in H^4(X,\mathbb{Q})$. Fix an integer $k$ such that $1\leq k \leq r$. Then 
 $$\overline{\Eff}^1\bigl(\Gr_X(k,E)\bigr) = \Bigl\{x\lambda_{E,k} + \pi_k^*\gamma \mid x\in \mathbb{R}_{\geq 0}, \gamma\in \overline{\Eff}^1(X)\Bigr\},$$
 where $\pi_k: \Gr_X(k,E) \longrightarrow X$ is the projection and $\lambda_{E,k} := c_1(\mathcal{O}_{\Gr_X(k,E)}(1)) - \frac{1}{\rk(\wedge^kE)}\pi_k^*c_1(\wedge^kE)$.
 
 In particular, if  $\overline{\Eff}^1(X)$ is a finite polyhedra generated by $\{E_1,E_2,\cdots,E_n\}$,
 then $$\overline{\Eff}^1\bigl(\Gr_X(k,E)\bigr) = \Bigl\{x_0\lambda_{E,k} + x_1\pi_k^*E_1+x_2\pi_k^*E_2+\cdots+x_n\pi_k^*E_n\mid x_i \in \mathbb{R}_{\geq 0}\Bigr\}.$$
 \end{thm}
 Miyaoka \cite{Miy87} showed that a vector bundle $E$ on a smooth curve $C$ is slope semistable if and only if $\Nef^1(\mathbb{P}_C(E)) = \overline{\Eff}^1(\mathbb{P}_C(E))$. This result is generalised in \cite{B-H-P} to Grassmann bundle $\Gr_C(k,E)$ over smooth projective curve $C$ defined over an algebraically closed field of positive characteristic, i.e. $E$ is strongly slope semistable on a smooth curve $C$ in characteristic $p>0$ if and only if $\Nef^1(\Gr_C(k,E)) = \overline{\Eff}^1(\Gr_C(k,E))$ for every $k$ satisfying $1\leq k \leq r$.
 However, examples are given in \cite{B-H-P} to show that this characterization of slope semistability in terms of equality of nef and pseudo-effective cones of Grassmann bundles $\Gr_X(k,E)$ over higher dimensional varieties $X$
is not true for vector bundles on higher dimensional projective varieties (see Section 5, \cite{B-H-P}).
 In this context, we prove the following equivalence:
 \begin{thm}
 Let $E$ be a semistable vector bundle of rank $r$ on a smooth irreducible projective surface $X$ such that $c_2\bigl(\End(E)\bigr) = 2rc_2(E)-(r-1)c^2_1(E) = 0 \in H^4(X,\mathbb{Q})$. Then the following are equivalent 
 
 \begin{enumerate}
 \item  $\overline{\Eff}^1(X) = \Nef^1(X)$.
 
 \item $c_1\bigl({\pi_k}_*\mathcal{O}_{\Gr_X(k,E)}(D)\bigr)\in \Nef^1(X)$ for every effective divisor $D$ in $\Gr_X(k,E)$. 

 \item  $\overline{\Eff}^1\bigl(\Gr_X(k,E)\bigr) = \Nef^1\bigl(\Gr_X(k,E)\bigr)$ 
 \end{enumerate}
\end{thm}
Next in Section 5, we give a description of the closed cone $\overline{\NE}\bigl(\Gr_X(k,E)\bigr)$ of curves in $\Gr_X(k,E)$  in terms of the closed cone of curves $\overline{\NE}\bigl(\Gr_C(k,E\vert_C)\bigr)$ for every irreducible curve $C$ in $X$. We compute the nef cone  $\Nef^1(\Gr_X(k,E))$ by applying duality in some special cases. In general, when the Picard number is at least 3, the nef cones might not be a finite polyhedron, and hence are not so easy to calculate. For example if $E$ is a rank 2 bundle obtained by the
Serre construction from the ideal sheaf of 10 very general points on $\mathbb{P}^2$, then the positivity
of $E$ is related to the Nagata conjecture. Thus one has to settle for special bundles $E$ on
special varieties $X$, even when dimension of $X$ is 2.  We mention few important findings here.

\begin{corl}
Let $X$ be a smooth complex projective surface with
\begin{center} 
$\overline{\NE}(X) = \Bigl\{ a_1[C_1] + a_2[C_2] +\cdots+a_n[C_n] \mid a_i \in \mathbb{R}_{\geq 0}\Bigr\}$
\end{center}
for  some irreducible curves $C_1,C_2,\cdots,C_n$  in $X$. If $E$ is a semistable vector bundle of rank $r\geq2$ on $X$ with $c_2\bigl(\End(E)\bigr) := 2rc_2(E)-(r-1)c_1^2(E) = 0 \in H^4(X,\mathbb{Q})$, then for any integer $k$ with $1\leq k \leq r$, we have
\begin{center}
 $\Nef^1\bigl(\Gr_X(k,E)\bigr)
   =\Bigl\{ y_0\xi+\pi_k^*\gamma \mid \gamma\in N^1(X)_{\mathbb{R}},  y_0\geq 0,    ky_0\mu(E\vert_{C_j}) + (\gamma\cdot C_j)\geq0 \text{ for all } 1\leq j\leq n\Bigr\}$,
\end{center}
where $\xi$  denotes the numerical equivalence class of 
the tautological bundle $\mathcal{O}_{\Gr(k,E)}(1)$.
\end{corl}
\begin{corl}
Let $X$ be a smooth complex projective surface with Picard number 1 and $L_X$ be an ample generator of the real N\'{e}ron Severi group $N^1(X)_{\mathbb{R}}$. 
 Let $E = M_1\oplus M_2\oplus\cdots\oplus M_r$ be a completely decomposable vector bundle of rank $r\geq2$ on $X$. 
  We fix an integer $k$ with $1\leq k \leq r$.
  Then
 \begin{center}
   $\Nef^1\bigl(\Gr_X(k,E)\bigr) = \Bigl\{ y_0\xi + y_1\pi^*L_X \mid y_0\geq 0, y_0\theta^{L_X}_{E,k}+y_1L_X^2 \geq 0\Bigr\}$
  \end{center}
 In particular, if $X = \mathbb{P}^2$, and $E = \mathcal{O}_{\mathbb{P}^2}(a_1)\oplus \mathcal{O}_{\mathbb{P}^2}(a_2) \oplus\cdots\oplus \mathcal{O}_{\mathbb{P}^2}(a_r)$
 a completely decomposable bundle on $\mathbb{P}^2$, then 
 \begin{align*}
 \Nef^1\bigl(\Gr_{\mathbb{P}^2}(k,E)\bigr) = \Bigl\{ y_0\xi + y_1\pi^*H \mid y_0\geq 0, y_0\theta_{E,k}^{H} + y_1\geq 0\Bigr\},
 \end{align*} 
 where $H$ is the numerical equivalence class of $\mathcal{O}_{\mathbb{P}^2}(1)$ in $N^1(\mathbb{P}^2)_{\mathbb{R}}$ (see subsection \ref{direct_sum} for the notation $\theta_{E,k}^{H}$).
 \end{corl}
We also compute the nef and pseudoeffective cones of the fiber products $\Gr_C(k,E)\times_C\Gr_C(k',E)$ over a smooth curve $C$. We give several examples and corollaries of our result.

Next in section 6, we prove the following:
\begin{thm}\label{thm1.5}
  Let $E$ be a vector bundle of rank $r$ on a smooth complex projective variety $X$. For every $k$ with $1\leq k\leq r$, consider the universal quotient bundle $Q_k$ on the grassmann bundle $\Gr_X(k,E)$ over $X$. Consider the following commutative diagram:
  \begin{center}
   \begin{tikzcd}
\mathbb{P}_Y(Q_k) \arrow[rd, "\phi_k"] \arrow[r, "\psi_k"] & Y=\Gr_X(k,E) \arrow[d,"\pi_k"]\\
& X
\end{tikzcd}
  \end{center}
Then the following are equivalent:
  \begin{enumerate}
   \item $\Phi_{E,k} := c_1(\mathcal{O}_{\mathbb{P}(Q_k)}(1)) - \frac{1}{r}\phi_k^*c_1(E)$ is nef for every $k$ with $1\leq k\leq r.$
   \item $E$ is semistable with $c_2\bigl(\End(E)\bigr)=0\in H^4(X,\mathbb{Q})$.
  \end{enumerate}
  \end{thm}
  Using Theorem \ref{thm1.5} we discuss about the nefness and ampleness of the universal quotient bundles. A particular application of Theorem \ref{thm1.5} is the following :
  \begin{corl}
 Let $E$ be a semistable  nef vector bundle of rank $r$ on a smooth complex projective variety $X$ of dimension $n$ with $c_2\bigl(\End(E)\bigr) = 0 \in H^4(X,\mathbb{Q})$. Then for every $k$ with $1\leq k \leq r$, the universal quotient bundle $Q_k$ is  nef on the grassmann bundle $\Gr_X(k,E)$.
 \end{corl}
 \section{Notation and Convention}
Throughout this article, all the algebraic varieties are assumed to be irreducible and reduced, and defined over the field of complex numbers $\mathbb{C}$.

Let $E$ be a vector bundle of rank $r$ on a projective variety $X$. For any positive integer $k$ with $1\leq k \leq r$, let  $\pi_k :\Gr(k,E)\longrightarrow X$ be the Grassmann bundle over $X$ associated to $E$ parametrizing all the quotients of dimension $k$ of the fibers of $E$.
In particular, we get the associated projective bundle $\mathbb{P}_X(E)$ when $k=1$, i.e. $\Gr(1,E) = \mathbb{P}_X(E)$.
We will simply write $\Gr(k,E)$ whenever the base space $X$ is clear from the context. The numerical class of the tautological line bundle $\mathcal{O}_{\Gr(k,E)}(1)$ will be denoted by $\xi$ unless otherwise specified.

The rank of a vector bundle $E$ will be denoted by $\rk(E)$. The numerical class of a divisor $D$ will be denoted by  $\bigl[D\bigr]$. The dual of a vector bundle $E$ will be denoted by $E^{\vee}$. In this article, semistability of a vector bundle means slope semistability unless otherwise specified.

The set of non-negative (resp. positive) real numbers will be denoted by $\mathbb{R}_{\geq 0}$ (resp. $\mathbb{R}_+$).
\section{Preliminaries}
 In this section we recall basic definition of nef and pseudo-effective cones and its various properties. We also recall the characterization of semistable bundles with vanishing discriminant over higher dimensional projective varieties which we will use to prove one of our main results. Finally, we recall the computation of nef and pseudoeffective cones of Grassmann bundles over smooth complex projective  curves.
 \subsection{Nef cone and Pseudo-effective cone}

 Let $X$ be a projective variety of dimension $n$ and $\mathcal{Z}_m(X)$ (respectively $\mathcal{Z}^m(X)$) denotes the free abelian group generated by $m$-dimensional (respectively $m$-codimensional) subvarieties on $X$. The Chow groups $A_m(X)$ are defined as the quotient of $\mathcal{Z}_m(X)$ modulo rational equivalence. When $X$ is a  smooth irreducible projective variety, we denote $A^m(X) := A_{n-m}(X)$.
 
 Two cycles $Z_1$, $Z_2 \in \mathcal{Z}^m(X)$ are said to be numerically equivalent, denoted by $Z_1\equiv Z_2$ if $Z_1\cdot \gamma =  Z_2\cdot\gamma $ 
  for all $\gamma \in \mathcal{Z}_m(X)$. The \it numerical groups \rm  $ N^m(X)_{\mathbb{R}}$ are defined as the quotient of 
  $\mathcal{Z}^m(X) \otimes \mathbb{R}$ modulo numerical equivalence. When $X$ is smooth, we define 
  $N_m(X)_{\mathbb{R}} := N^{n-m}(X)_{\mathbb{R}}$ for all $m$.
  
  Let 
  \begin{center}
  $\Div^0(X) := \Bigl\{ D \in \Div(X) \mid D\cdot C = 0 $ for all curves $C$ in $X \Bigr\} \subseteq \Div(X)$.
  \end{center}
 be the subgroup of $\Div(X)$ consisting of numerically trivial divisors. The quotient $\Div(X)/\Div^0(X)$ is called the N\'{e}ron Severi group of $X$, and is denoted by $N^1(X)_{\mathbb{Z}}$.
   The N\'{e}ron Severi group  $N^1(X)_{\mathbb{Z}}$ is a free abelian group of finite rank.
 Its rank, denoted by $\rho(X)$ is called the Picard number of $X$. In particular, $N^1(X)_{\mathbb{R}}$ is called the real N\'{e}ron 
 Severi group and $N^1(X)_{\mathbb{R}}  := N^1(X)_{\mathbb{Z}} \otimes \mathbb{R} := \bigl(\Div(X)/\Div^0(X)\bigr) \otimes \mathbb{R}$.
 For $X$ smooth, the intersection product induces a perfect pairing 
 \begin{align*}
N^m(X)_{\mathbb{R}} \times N_m(X)_{\mathbb{R}} \longrightarrow \mathbb{R}
\end{align*}
which implies $N^m(X)_{\mathbb{R}} \simeq (N_m(X)_{\mathbb{R}})^\vee$ for every $m$ satisfying $0\leq m\leq n$. The direct sum 
  $N(X) := \bigoplus\limits_{m=0}^{n}N^m(X)_{\mathbb{R}}$ is a graded $\mathbb{R}$-algebra with multiplication induced by the intersection form.
  
  The convex cone generated by the set of all effective $m$-cycles in $N_m(X)_\mathbb{R}$ is denoted by $\Eff_m(X)$ and its closure $\overline{\Eff}_m(X)$ is called the \it pseudo-effective cone \rm  of $m$-cycles in $X$. For any $0 \leq m \leq n$,
$\overline{\Eff}^m(X) := \overline{\Eff}_{n-m}(X)$. We denote the closed cone of curves $\overline{\Eff}_1(X)$ by $\overline{\NE}(X)$. 

The \it nef cone \rm  are defined as follows :
\begin{align*}
 \Nef^m(X) := \Bigl\{ \alpha \in N^m(X) \mid \alpha \cdot \beta \geq 0 \hspace{2mm} \forall \beta \in \overline{\Eff}_m(X)\Bigr\}.  
\end{align*}
In particular, $\Nef^1(X)$ is dual to the closed cone of curves $\overline{\NE}(X)$.

If $p : U\longrightarrow  W$ is a family of $m$-cycles on $X$ and $U$ is irreducible and the second projection $s : U \longrightarrow X$ is
dominant, we say that $p$ is \it strongly movable\rm, and that $[p]$ is a \it strongly movable class \rm represented by the family
$p$. If $U$ is reducible, but every irreducible component still dominates $X$, we say that $p$ is strictly movable.
The closure of the cone in $N_m(X)$ generated by strongly (or strictly) movable classes is called the \it movable cone \rm $\overline{\Mov}_m(X)$.

In particluar, an irreducible curve $C$ in $X$ is called \it movable \rm if there exists an algebraic family of irreducible curves $\{C_t\}_{t\in T}$ such 
that $C = C_{t_0}$ for some $t_0 \in T$ and $\bigcup_{t \in T} C_t \subset X$ is dense in $X$.

A class $\gamma \in N_1(X)_{\mathbb{R}}$ is called movable if there exists a movable curve $C$
such that $\gamma = [C]$ in $N_1(X)_{\mathbb{R}}$. The closure of the cone generated by movable classes in
$N_1(X)_{\mathbb{R}}$, denoted by $\overline{\Mov}_1(X)$ is called the \it movable cone\rm. By \cite{BDPP13} $\overline{\Mov}_1(X)$ is 
the dual cone to $\overline{\Eff}^1(X)$. Also, $\Nef^1(X) \subseteq \overline{\Eff}^1(X)$ always (see Ch 2 \cite{L1}). We refer the 
reader to \cite{L1},\cite{L2} for more details about these cones.

\subsection{Semistability of Vector bundles}
Let $X$ be a smooth complex projective variety of dimension $n$ with a fixed ample line bundle $H$ on it.
For a torsion-free coherent sheaf $\mathcal{G}$ of rank $r$ on $X$, the degree with respect to $H$, denoted by $\deg_H(\mathcal{G})$ is defined as follows :
$$\deg_H(\mathcal{G}) =c_1(\mathcal{G})\cdot H^{n-1},$$
and the $H$-slope of $\mathcal{G}$ is defined as 
\begin{align*}
\mu_H(\mathcal{G}) := \frac{c_1(\mathcal{G})\cdot H^{n-1}}{r} \in \mathbb{Q}.
\end{align*}
A vector bundle  $E$ on $X$ is said to be $H$-semistable if $\mu_H(\mathcal{G}) \leq \mu_H(E)$ for 
all subsheaves $\mathcal{G}$ of $E$.
A vector bundle $E$ on $X$ is called $H$-unstable if it is not $H$-semistable. For every vector bundle $E$ on $X$, there is a unique filtration
\begin{align*}
 0 = E_0 \subsetneq E_{1} \subsetneq E_{2} \subsetneq\cdots\subsetneq E_{l-1} \subsetneq E_l = E
\end{align*}
of subbundles of $E$, called the Harder-Narasimhan filtration of $E$, such that $E_i/E_{i-1}$ is $H$-semistable torsion-free sheaf
for each $i \in \{ 1,2,\cdots,l\}$
and $\mu_H\bigl(E_{1}/E_{0}\bigr) > \mu_H\bigl(E_{2}/E_{1}\bigr) >\cdots> \mu_H\bigl(E_{l}/E_{l-1}\bigr)$.
We define $Q_1:= E_{l}/E_{l-1}$ and $\mu_{\min}(E) := \mu_H(Q_1) = \mu_H\bigl(E_{l}/E_{l-1}\bigr)$ and $\mu_{\max}(E) := \mu_H\bigl(E_{l}/E_{l-1}\bigr).$

For a vector bundle $E$ of rank $r$ over $X$, the element $c_2\bigl(\End(E)\bigr) \in H^4(X,\mathbb{Q})$ is called the \it discriminant \rm of $E$.  
We recall the following result from \cite{N99}  or  [Theorem 1.2,\cite{B-B}] about semistable vector bundles with vanishing discriminant.
\begin{thm}\label{thm3.1}
 Let $E$ be a vector bundle of rank $r$ on a smooth complex projective variety $X$. Let $\pi : \mathbb{P}(E)\longrightarrow X$ be the projection map. Then the following are equivalent
 
 \rm(1) \it $E$ is semistable and $c_2\bigl(\End(E)\bigr) = 0$.
 
 \rm(2) \it $\lambda_E := c_1(\mathcal{O}_{\mathbb{P}(E)}(1)) - \frac{1}{r}\pi^*c_1(E) \in \Nef^1\bigl(\mathbb{P}(E)\bigr)$.
 
 \rm(3) \it For every pair of the form $(\phi,C)$, where $C$ is a smooth projective curve and $\phi : C \longrightarrow X$ is a non-constant morphism, $\phi^*(E)$ is semistable.
\end{thm}
 Since nefness of a line bundle does not depend on the fixed polarization (i.e. fixed ample line bundle $H$) on $X$, the Theorem \ref{thm3.1} implies that the semistability of a vector bundle $E$ with $c_2\bigl(\End(E)\bigr) = 0$ is independent of the fixed polarization $H$. We will not mention about the polarization $H$ from now on whenever we speak of semistability of such bundles $E$ with $c_2\bigl(\End(E)\bigr)=0$.
We have the following lemmata as easy applications of Theorem \ref{thm3.1}.
\begin{lem}\label{lem3.2}
 Let $\psi : X \longrightarrow Y$ be a morphism between two smooth complex projective varieties and $E$ is a semistable bundle on $Y$ with $c_2\bigl(\End(E)\bigr) = 0$. Then the pullback bundle $\psi^*(E)$ is also semistable with $c_2\bigl(\End(\psi^*E)\bigr) = 0$.
 \begin{proof}
  Let $\phi : C \longrightarrow X$ be a non-constant morphism from a smooth curve $C$ to $X$. If the image of $\phi$ is contained in any fibre of $\psi$, then the pullback bundle $\phi^*\psi^*(E)$ is trivial, and hence semistable. Now let us assume that the image is not contained in any fibre of $\phi$. As $E$ is semistable bundle on $Y$ with $c_2\bigl(\End(E)\bigr) = 0$, by the previous Theorem \ref{thm3.1} the pullback bundle  $(\psi\circ\phi)^*E = \phi^*(\psi^*E)$ under the non-constant morphism $\psi\circ\phi$  is semistable on $C$. Hence $\psi^*E$ is semistable bundle on $X$ with $c_2\bigl(\End(\psi^*E)\bigr) = 0$.
 \end{proof}
\end{lem}
\begin{lem}\label{lem3.3}
  Let  $E$ be a semistable vector bundle of rank $r$ on a smooth complex projective variety $X$ with $c_2\bigl(\End(E)\bigr) = 0$. Then for any positive integers $m$ and $k$, the vector bundle $\bigl(\wedge^kE\bigr)^{\otimes m}$ is semistable with $c_2\bigl(\End\bigl(\wedge^kE\bigr)^{\otimes m}\bigr)=0$.
\begin{proof}
  Note that for any smooth complex projective curve $C$ and any non-constant map $\phi : C \longrightarrow X$, the pullback bundle $\phi^*\bigl(\wedge^k(E)\bigr)^
 {\otimes m} = \bigl(\wedge^k(\phi^*(E))\bigr)^
 {\otimes m}$ is semistable bundle on $C$. Therefore, by Theorem \ref{thm3.1}, we have that the vector bundle $\bigl(\wedge^kE\bigr)^{\otimes m}$ is semistable with $c_2\bigl(\End\bigl(\wedge^kE\bigr)^{\otimes m}\bigr)=0$.
\end{proof}
\end{lem}
\begin{lem}\label{lem3.4}
 Let $E$ be a semistable vector bundle of rank $r$ on a smooth complex projective variety $X$ of dimension $d$ with $c_2\bigl(\End(E)
 \bigr) = 0$. Fix two integers $m$ and $k$ such that $m\geq 0$ and $1\leq k \leq r$. Consider the projection $\pi : \Gr(k,E) \longrightarrow X$. Then the vector bundle $\pi_*(\mathcal{O}_{\Gr(k,E)}(m))$ is semistable with $c_2\bigl(\End(\pi_*(\mathcal{O}_{\Gr(k,E)}(m)))\bigr) = 0$.
 \begin{proof}
  We note that $E_{k,m} := \pi_*(\mathcal{O}_{\Gr(k,E)}(m))$ is the vector bundle associated to $E$ by the Weyl module with highest weight $m\omega_k$ with $\omega_k$ as the fundamental weight corresponding to the $k$-th exterior power presentation. Hence $E_{k,m}$ is a direct summand of $\bigl(\wedge^kE\bigr)^{\otimes m}$.
  
  Hence for any ample bundle $H$ on $X$, we have $$\mu_H(E_{k,m}) = \mu_H\bigl(\bigl(\wedge^kE\bigr)^{\otimes m}\bigr).$$
  
  This implies $\mu_H\bigl(E_{k,m}\otimes ((\wedge^kE)^{\otimes m})^{\vee}\bigr) = \mu_H\bigl(E_{k,m}) - \mu_H\bigl(\bigl(\wedge^kE\bigr)^{\otimes m}\bigr) =0$.
  
  Hence we have
  $c_1\bigl(E_{k,m}\otimes ((\wedge^kE)^{\otimes m})^{\vee}\bigr)\cdot H^{d-1} = 0$ for any ample line bundle $H$ on $X$.
  
  This shows that $c_1\bigl(E_{k,m}\otimes ((\wedge^kE)^{\otimes m})^{\vee}\bigr) \equiv 0$ in $N^1(X)_{\mathbb{R}}$.
  
  Let  $\phi : C \longrightarrow X$ be a non-constant morphism from a smooth curve $C$ to $X$. 
  
  Now consider  the following exact sequence $$0\longrightarrow K\longrightarrow (\wedge^kE)^{\otimes m} \longrightarrow E_{k,m}\longrightarrow 0$$
  Then we have $$0\longrightarrow \phi^*(K)\longrightarrow \phi^*(\wedge^kE)^{\otimes m} \longrightarrow \phi^*(E_{k,m})\longrightarrow 0.$$
  
  By the above observation we have $\mu\bigl(\phi^*(K)\bigr) = \mu\bigl(\phi^*(\wedge^kE)^{\otimes m}\bigr) = \mu\bigl( \phi^*(E_{k,m})\bigr)$. 
  
  This shows that $\phi^*(E_{k,m})$ is a semistable bundle on $C$.  Therefore  we conclude that $E_{k,m}$ is a semistable bundle with $c_2\bigl(\End(E_{k,m})\bigr) = 0$.
 \end{proof}
\end{lem}
\subsection{Positive cones of Grassmann Bundles over curves}
Let $E$ be a vector bundle of rank $r$ on a smooth complex projective curve $C$. Let $$0=E_0\subsetneq E_1\subsetneq E_2\subsetneq \cdots\subsetneq E_{l-1} \subsetneq E_l =E$$ be the Harder-Narasimhan filtration of $E$. For an integer $k$ with $1\leq k \leq r$, the nef cone $\Nef^1(\Gr_C(k,E))$ and pseudoeffective cone $\overline{\Eff}^1(\Gr_C(k,E))$ has been computed in \cite{B-P} and \cite{B-H-P} respectively. We recall these results here.

Define $$\theta_{E,k} := (k-\rk(E/E_t))\mu(E_t/E_{t-1})+\deg(E/E_t),$$ where either $t=l$ or $t$ is the smallest integer with $\sum\limits_{i=t+1}^d\rk(E_i/E_{i-1}) = \rk(E/E_t) < k$.

We then have $$\Nef^1(\Gr(k,E)) = \Bigl\{a(\xi - \theta_{E,k}f) +bf\mid a,b\in\mathbb{R}_{\geq  0}\Bigr\},$$ where $\xi$ is the numerical class of the tautological line bundle $\mathcal{O}_{\Gr(k,E)}(1)$ and $f$ denotes the fiber of the projection map $\pi : \Gr(k,E)\longrightarrow C$.  Note that  $\theta_{E,k} =k\mu(E)$ when $E$ is semistable.

Define $$\zeta_{E,k} : =(k-\rk(E_t))\mu(E_{t+1}/E_t)+\deg(E_t),$$ where $t$ is the unique smallest integer such that $\sum\limits_{i=1}^{t+1}\rk(E_i/E_{i-1})>k$, So if $\rk(E_1)>k$, then $t=0$; otherwise $t$ is the largest integer such that $\sum\limits_{i=1}^t\rk(E_i/E_{i-1}) \leq k$. 

We then have $$\overline{\Eff}^1(\Gr(k,E))=\Bigl\{a(\xi-\zeta_{E,k}f)+bf\mid a,b\in \mathbb{R}_{\geq0}\Bigr\}.$$

Note that if $E$ is semistable, then $t=0$. In that case, $\Nef^1(\Gr(k,E))= \overline{\Eff}^1(\Gr(k,E))$ for every $k$ with $1\leq k\leq r$.
\section{pseudoeffective cone of grassmann bundle}
\begin{thm}\label{thm4.1}
 Let $E$ be a semistable vector bundle of rank $r$ on a smooth irreducible complex projective variety $X$ such that $c_2\bigl(\End(E)\bigr) = 2rc_2(E)-(r-1)c^2_1(E) = 0 \in H^4(X,\mathbb{Q})$. Fix an integer $k$ such that $1\leq k \leq r$. Then 
 $$\overline{\Eff}^1\bigl(\Gr(k,E)\bigr) = \Bigl\{x\lambda_{E,k} + \pi_k^*\gamma \mid x\in \mathbb{R}_{\geq 0}, \gamma\in \overline{\Eff}^1(X)\Bigr\},$$ where $\pi_k:\Gr(k,E)\longrightarrow X$ is the projection, and 
 $\lambda_{E,k} = c_1\bigl(\mathcal{O}_{\Gr(k,E)}(1)\bigr) - \frac{1}{\rk(\wedge^k(E)}\pi_k^*c_1\bigl(\wedge^kE)\bigr).$
 
 In particular, if $\overline{\Eff}^1(X)$ is a finite polyhedra generated by $\{E_1,E_2,\cdots,E_n\}$.
 Then $$\overline{\Eff}^1\bigl(\Gr(k,E)\bigr) = \Bigl\{x_0\lambda_{E,k} + x_1\pi_k^*E_1+x_2\pi_k^*E_2+\cdots+x_n\pi^*_kE_n\mid x_i \in \mathbb{R}_{\geq 0}\Bigr\}.$$
\end{thm}
\begin{proof}
Let $D$ be an effective divisor on $\Gr(k,E)$ such that 
$\mathcal{O}_{\Gr(k,E)}(D) \cong \mathcal{O}_{\Gr(k,E)}(m) \otimes \pi_k^*(\mathcal{L})$ for some integer $m$ and a line bundle $\mathcal{L} \in \Pic(X)$. Then 
 \begin{align*}
  H^0\bigl(\Gr(k,E), \mathcal{O}_{\Gr(k,E)}(m) \otimes \pi_k^*(\mathcal{L})\bigr) = H^0\bigl(X, E_{k,m} \otimes \mathcal{L}\bigr) \neq 0 ,
 \end{align*}
 which implies $m \geq 0$ (Here $E_{k,m} := {\pi_k}_*(\mathcal{O}_{\Gr(k,E)}(m))$ is the vector bundle associated to $E$ by the Weyl module with highest weight $m\omega_k$ with $\omega_k$ as the fundamental weight corresponding to the $k$-th exterior power presentation). 
 
Let $\gamma = [C] $ be a movable class in $N_1(X)_{\mathbb{R}}$. Then $C$ belongs to an algebraic family of curves
$\{C_t\}_{t\in T}$ such that  $\bigcup_{t \in T} C_t$ covers a dense subset of $X$. So we can find a curve
$C_{t_1}$ in this family such that  
\begin{align*}
 H^0\bigl( C_{t_1}, E_{k,m}\vert_{C_{t_1}} \otimes \mathcal{L}\vert_{C_{t_1}}\bigr) \neq 0.
\end{align*}

Let $\eta_{t_1} :\tilde{C_{t_1}} \longrightarrow C_{t_1}$ be the normalization of the curve $C_{t_1}$ and we call $\phi_{t_1} :=  i\circ \eta_{t_1} $ where $i : C_{t_1} \hookrightarrow X$ is the inclusion.

As $E$ is a semistable bundle on $X$ with $c_2\bigl(\End(E)\bigr) = 0$, by Lemma \ref{lem3.4} $E_{k,m} \otimes \mathcal{L}$ is also semistable on $X$ and $c_2\bigl(\End(E_{k,m} \otimes \mathcal{L})\bigr) = 0$. Therefore  we have $\phi_{t_1}^*\bigl(E_{k,m} \otimes \mathcal{L}\bigr)$ is also semistable on $\tilde{C_{t_1}}$. 

Since $\eta_{t_1}$ is a surjective map, we have
\begin{align*}
 H^0\bigl(C_{t_1}, \eta_{t_1}^*(E_{k,m}\vert_{C_{t_1}} \otimes \mathcal{L}\vert_{C_{t_1}})\bigr)  = H^0\bigl(C_{t_1}, \phi_{t_1}^*(E_{k,m} \otimes \mathcal{L})\bigr)\neq 0.
\end{align*}
This implies that $c_1\bigl(\phi_{t_1}^*\bigl(E_{k,m} \otimes \mathcal{L}\bigr)\bigr) \geq 0$ , and hence

$c_1\bigl(E_{k,m} \otimes \mathcal{L}\bigr)\cdot C_{t_1} = \bigl\{ c_1\bigl(E_{k,m} \otimes \mathcal{L}\bigr) \cdot \gamma \bigr\} \geq 0$ for a movable class $\gamma \in N_1(X)_{\mathbb{R}}$.

Using the duality property of movable cone $\overline{\Mov}_1(X)$ we conclude that
\begin{center}
$ c_1\bigl(E_{k,m} \otimes \mathcal{L}\bigr) = c_1\bigl({\pi_k}_*\mathcal{O}_{\Gr(k,E)}(D)\bigr) \in \overline{\Eff}^1(X)$.
\end{center}
Now 
 $c_1\bigl({\pi_k}_*\mathcal{O}_{\Gr(k,E)}(D)\bigr)$
 = $c_1\bigl(E_{k,m}\otimes \mathcal{L}\bigr) =  c_1(E_{k,m}) + \rk(E_{k,m}) c_1(\mathcal{L})  \in \overline{\Eff}^1(X)$.
 
 Note that as in Lemma \ref{lem3.4}, we have $c_1\bigr(E_{k,m}\otimes {((\wedge^kE)^{\otimes m}})^{\vee}\bigl) \equiv 0$, so that $$\frac{c_1(E_{k,m})}{\rk(E_{k,m})}+ c_1(\mathcal{L}) = \frac{m}{\rk(\wedge^kE)}c_1\bigl(\wedge^kE\bigr) + c_1(\mathcal{L}) \in \overline{\Eff}^1(X).$$

 Since $E$ is a semistable bundle on $X$ with $c_2\bigl(\End(E)\bigr) = 0$,  we also have $\wedge^kE$ is semistable with $c_2\bigl(\End(\wedge^k(E))\bigr) = 0$. Hence we have $$c_1\bigl(\mathcal{O}_{\mathbb{P}(\wedge^kE)}(1)\bigr) - \frac{1}{\rk(\wedge^k(E))}p_k^*\bigl(c_1(\wedge^kE)\bigr) \in \Nef^1\bigl(\mathbb{P}(\wedge^kE)\bigr) \subseteq \overline{\Eff}^1\bigl(\mathbb{P}(\wedge^kE)\bigr),$$
 where $p_k : \mathbb{P}_X(\wedge^k E) \longrightarrow X$ is the projection map.
 
 This implies $$\lambda_{E,k} = c_1\bigl(\mathcal{O}_{\Gr(k,E)}(1)\bigr) - \frac{1}{\rk(\wedge^k(E)}\pi_k^*c_1\bigl(\wedge^kE)\bigr) \in \Nef^1\bigl(\Gr(k,E)\bigr) \subseteq \overline{\Eff}^1\bigl(\Gr(k,E)\bigr).$$
Hence,
$\mathcal{O}_{\Gr(k,E)}(m) \otimes \pi_k^*(\mathcal{L})$
$\equiv m\bigl(c_1(\mathcal{O}_{\Gr(k,E)}(1))\bigr) + \pi_k^*c_1(\mathcal{L})$

$\equiv m\Bigl\{ c_1(\mathcal{O}_{\Gr(k,E)}(1)) - \frac{1}{\rk(\wedge^k(E)}\pi_k^*c_1\bigl(\wedge^kE\bigr) \Bigr\} + \frac{m}{\rk(\wedge^k(E))}\bigl(\pi_k^*c_1(\wedge^kE) + \pi_k^*c_1(\mathcal{L})\bigr)$

$\equiv m\Bigl\{ c_1(\mathcal{O}_{\Gr(k,E)}(1)) - \frac{1}{\rk(\wedge^k(E)}\pi_k^*c_1\bigl(\wedge^kE\bigr)\Bigr\} + \pi_k^*\bigl(\frac{m}{\rk(\wedge^kE)}c_1(\
\wedge^kE) + c_1(\mathcal{L})\bigr) \in \overline{\Eff}^1\bigl(\Gr(k,E)\bigr)$.

This shows that 
$$ \overline{\Eff}^1\bigl(\Gr(k,E)\bigr) = \Bigl\{x\lambda_{E,k} + \pi_k^*\gamma \mid x\in \mathbb{R}_{\geq 0}, \gamma\in \overline{\Eff}^1(X)\Bigr\}. $$
Moreover, if $\overline{\Eff}^1(X)$ is a finite polyhedra generated by $E_1,E_2,\cdots, E_n$, then 
$$\frac{m}{\rk(\wedge^kE)}c_1(\wedge^kE)+c_1(\mathcal{L}) = y_1E_1+y_2E_2+\cdots+y_nE_n$$ for some $y_1,y_2,\cdots,y_n\in \mathbb{R}_{\geq 0}$, say.

Thus 
$\mathcal{O}_{\Gr(k,E)}(m) \otimes \pi_k^*(\mathcal{L})$
$\equiv m\bigl(c_1(\mathcal{O}_{\Gr(k,E)}(1))\bigr) + \pi_k^*c_1(\mathcal{L})$

$\equiv m\bigl\{ c_1(\mathcal{O}_{\Gr(k,E)}(1)) - \frac{1}{\rk(\wedge^kE)}\pi^*c_1(E)\bigr\} + y_1\pi_k^*E_1+y_2\pi_k^*E_2+\cdots+y_n\pi_k^*E_n $
\vspace{2mm}

Therefore, $\Eff^1\bigl(\Gr(k,E)\bigr) \subseteq \Bigl\{x_0\lambda_{E,k} + x_1\pi_k^*E_1+x_2\pi_k^*E_2+\cdots+x_n\pi_k^*E_n\mid x_i \in \mathbb{R}_{\geq 0}\Bigr\} \subseteq \overline{\Eff}^1\bigl(\Gr(k,E)\bigr)$. 
\vspace{2mm}

Taking closure, we get $$\overline{\Eff}^1\bigl(\Gr(k,E)\bigr) = \Bigl\{x_0\lambda_{E,k} + x_1\pi_k^*E_1+x_2\pi_k^*E_2+\cdots+x_n\pi_k^*E_n\mid x_i \in \mathbb{R}_{\geq 0}\Bigr\}.$$
\end{proof}
We note that for a vector $E$ on a smooth irreducible projective complex variety $X$, the equality
$\Nef^1\bigl(\Gr(k,E)\bigr) = \overline{\Eff}^1\bigl(\Gr(k,E)\bigr)$ always implies the equality $\Nef^1(X) = \overline{\Eff}^1(X)$. However, 
the converse is not true in general (see section 5, \cite{B-H-P}). Next we prove that 
$\Nef^1\bigl(\Gr(k,E)\bigr) = \overline{\Eff}^1\bigl(\Gr(k,E)\bigr)$ if and only if $\Nef^1(X) = \overline{\Eff}^1(X)$ under the assumption that
$E$ is a semistable bundle on $X$ with $c_2\bigl(\End(E)\bigr) = 0.$ 
\begin{thm}
 Let $E$ be a semistable vector bundle of rank $r$ on a smooth irreducible projective surface $X$ such that $c_2\bigl(\End(E)\bigr) = 2rc_2(E)-(r-1)c^2_1(E) = 0 \in H^4(X,\mathbb{Q})$. Then the following are equivalent 
 
 \begin{enumerate}
 \item  $\overline{\Eff}^1(X) = \Nef^1(X)$.
 
 \item $c_1\bigl({\pi_k}_*\mathcal{O}_{\Gr(k,E)}(D)\bigr)\in \Nef^1(X)$ for every effective divisor $D$ in $\Gr(k,E)$. 

 \item  $\overline{\Eff}^1\bigl(\Gr(k,E)\bigr) = \Nef^1\bigl(\Gr(k,E)\bigr)$ 
 \end{enumerate}
\end{thm}
\begin{proof}
 (1)$\implies$(2) If $D$ is an effective divisor in  $\Gr(k,E)$, then procedding as in Theorem \ref{thm4.1}, we have $$c_1({\pi_k}_*\mathcal{O}_{\Gr(k,E)}(D)) \in \overline{\Eff}^1(X) = \Nef^1(X).$$
 
 (2)$\implies$(3) Suppose $D$ is an effective divisor in $\Gr(k,E)$ such that $$\mathcal{O}_{\Gr(k,E)}(D) \cong \mathcal{O}_{\Gr(k,E)}(m) \otimes \pi_k^*(\mathcal{L})$$ for some integer $m$ and a line bundle $\mathcal{L} \in \Pic(X)$. Then 
 \begin{align*}
  H^0\bigl(\Gr(k,E), \mathcal{O}_{\Gr(k,E)}(m) \otimes \pi_k^*(\mathcal{L})\bigr) = H^0\bigl(X, E_{k,m} \otimes \mathcal{L}\bigr) \neq 0 ,
 \end{align*}
 which implies $m \geq 0$. 
 
 Also by given hypothesis $c_1\bigl({\pi_k}_*\mathcal{O}_{\Gr(k,E)}(D)\bigr) = c_1(E_{k,m}\otimes \mathcal{L}) \in \Nef^1(X).$
 
 This implies $\frac{m}{\rk(\wedge^kE)}c_1(
\wedge^kE) + c_1(\mathcal{L}) \in \Nef^1(X)$.
 
 Hence,
$\mathcal{O}_{\Gr(k,E)}(m) \otimes \pi_k^*(\mathcal{L})$
$\equiv m\bigl(c_1(\mathcal{O}_{\Gr(k,E)}(1))\bigr) + \pi_k^*c_1(\mathcal{L})$

$\equiv m\Bigl\{ c_1(\mathcal{O}_{\Gr(k,E)}(1)) - \frac{1}{\rk(\wedge^k(E)}\pi_k^*c_1\bigl(\wedge^kE\bigr) \Bigr\} + \frac{m}{\rk(\wedge^k(E))}\bigl(\pi_k^*c_1(\wedge^kE) + \pi_k^*c_1(\mathcal{L})\bigr)$

$\equiv m\Bigl\{ c_1(\mathcal{O}_{\Gr(k,E)}(1)) - \frac{1}{\rk(\wedge^k(E)}\pi_k^*c_1\bigl(\wedge^kE\bigr)\Bigr\} + \pi_k^*\bigl(\frac{m}{\rk(\wedge^kE)}c_1(\
\wedge^kE) + c_1(\mathcal{L})\bigr) \in \Nef^1\bigl(\Gr(k,E)\bigr)$.
 
 (3)$\implies$(1) 
We claim that $\Nef^1(X) = \overline{\Eff}^1(X)$. If not then there is an effective divisor $D$ in $X$ which is not nef. Therefore the pullback $\pi_k^*(D)$ is also effective. Since $\pi_k$ is a surjective proper morphism, $\pi_k^*(D)$ can never be nef, which contradicts that $\Nef^1\bigl(\Gr(k,E)\bigr) = \overline{\Eff}^1\bigl(\Gr(k,E)\bigr).$ Hence we are done.
\end{proof}
\begin{corl}\label{corl4.3}
Let $X$ be a smooth complex projective variety of dimension $n$ with Picard number $\rho(X) = 1$ and $E$ be a semistable vector bundle of rank $r$ on $X$ with $c_2\bigl(\End(E)\bigr) = 0$. Fix an integer $k$ with $1\leq k \leq r$. Then
\begin{center}
$\Nef^1\bigl(\Gr(k,E)\bigr) = \overline{\Eff}^1\bigl(\Gr(k,E)\bigr) = \Bigl\{ x_0\lambda_{E,k} + x_1(\pi_k^*L_X) \mid x_0,x_1 \in \mathbb{R}_{\geq 0}\Bigr\}$,
\end{center}
 where $L_X$ is the numerical class of an ample generator of the Ne\'{r}on-Severi group $N^1(X)_{\mathbb{Z}}\cong \mathbb{Z}$.
 \begin{proof}
 Let $D$ be an effective divisor on $X$ and $D\equiv mL_X$ for some $m\in \mathbb{Z}$. Then $L_X^{n-1}\cdot D = mL_X^{n} \geq 0$. As $L_X$ is ample, by Nakai criterion for ampleness we have $L_X^{n} > 0$. This shows that $m \geq 0$ and $D \equiv mL_X \in \Nef^1(X)$. 
So $\Nef^1(X) = \overline{\Eff}^1(X)$ and hence the result follows from  Theorem \ref{thm4.1}.
\end{proof}
\end{corl}
\begin{xrem}\label{xrem1}
\rm In  [Corollary 4.3,\cite{B-H-P}], it is proved that a vector bundle $E$ of rank $r$ on a smooth complex projective curve $C$ is semistable if and only if for every $k$ with $1\leq k \leq r$, the following equality holds:
\begin{align*}
\Nef^1\bigl(\Gr(k,E)\bigr) = \overline{\Eff}^1\bigl(\Gr(k,E)\bigr) = \Bigl\{ x\lambda_{E,k} + y\pi_k^*L_C \mid x,y \in \mathbb{R}_{\geq 0}\Bigr\} = \Bigl\{ x\theta_{E,k} + yf\mid x,y \in \mathbb{R}_{\geq 0}\Bigr\},
\end{align*}
 where $L_C$ is the ample generator of $N^1(C)_{\mathbb{Z}}$ and $f$ is the fiber of the projection map $\pi_k : \Gr(k,E)\longrightarrow C$. The smooth curve $C$ has Picard number 1 and any second Chern class vanishes on $C$. So the Corollary \ref{corl4.3} can be thought as a partial generalization of  the above mentioned result in \cite{B-H-P} for vector bundles over higher dimensional varieties. 
\end{xrem}
\begin{exm}\label{exm4.4}
\rm Let $X$ be a smooth projective variety with Picard number $\rho(X) = 1$, and $E = \mathcal{L}_1\oplus \mathcal{L}_2 \oplus\cdots\oplus \mathcal{L}_r$ be a rank $r$ bundle on $X$ such that $\mathcal{L}_1\cdot L_X = \mathcal{L}_2\cdot L_X = \cdots = \mathcal{L}_r\cdot L_X$ (Here $L_X$ is the ample generator for $N^1(X)_{\mathbb{Z}}$). Then $E$ is semistable with $c_2\bigl(\End(E)\bigr) = 4c_2(E) - c_1^2(E) = 0$. Therefore, by Corollary \ref{corl4.3}  for every $k$ with $1\leq k \leq r$, we have $$\Nef^1\bigl(\Gr(k,E)\bigr) = \overline{\Eff}^1\bigl(\Gr(k,E)\bigr)= \Bigl\{x_0\lambda_{E,k}+x_1(\pi_k^*L_X)\mid x_0,x_1\in \mathbb{R}_{\geq 0}\Bigr\}.$$
\end{exm}
 \begin{exm}\label{exm4.5}
  \rm Let $G$ be a connected algebraic group acting transitively on a complex projective variety $X$. Then every effective divisor on $X$ is nef, i.e. $\Nef^1(X) = \overline{\Eff}^1(X)$. This is because any irreducible curve $C\subset X$ meets the translate $gD$ of an effective divisor $D$ properly for a general element $g \in G$. Since $G$ is connected, $gD \equiv D$. Therefore $D\cdot C = gD \cdot C \geq 0$, and hence $D$ is nef.  Examples of such homogeneous varieties $X$ include smooth abelian varieties, flag manifolds etc. So for any semistable vector bundle $E$ of rank $r$ on such a homogeneous space $X$ with $c_2\bigl(\End(E)\bigr) = 0$, using Theorem \ref{thm4.1} we have $\Nef^1\bigl(\Gr(k,E)\bigr) = \overline{\Eff}^1\bigl(\Gr(k,E)\bigr)$ for every integer $k$ with $1
  \leq k \leq r$. 
  
  For example, let $C$ be a general elliptic curve and $X = C\times C$ be the self product. Then $X$ is an abelian surface and $\Nef^1(X)$ is a non-polyhedral cone (see Lemma 1.5.4 in \cite{L1}). Let $p_i : X \longrightarrow C$ be the projection maps. For any semistable vector bundle $F$ of rank $r$ on $C$, the pullback bundle $E:= p_i^*(F)$ is a semistable bundle with $c_2\bigl(\End(E)\bigr) = 0$. Hence $\Nef^1\bigl(\Gr(k,E)\bigr) = \overline{\Eff}^1\bigl(\Gr(k,E)\bigr)$ for every integer $k$ with $1\leq k \leq r$.
 
\rm The above example can be extended as follows. Let $B$ be any smooth curve and $C$ be an elliptic curve. Then the product of $C$ with the Jacobian variety of $B$ i.e.,  $C \times J(B)$ is an abelian variety. For any semistable vector bundle $F$ of rank $r$ on $C$, the pullback bundle $E := p_1^*(F)$  under the 1st projection $p_1$, is a semistable bundle with $c_2\bigl(\End(E)\bigr) = 0$. Hence $\Nef^1\bigl(\Gr(k,E)\bigr) = \overline{\Eff}^1\bigl(\Gr(k,E)\bigr)$ for every integer $k$ with $1\leq k\leq r.$
 \end{exm}

 A vector bundle $E$ on an abelian variety $X$ is called weakly-translation invariant (semi-homogeneous in the sense of Mukai) if 
 for every closed point $x \in X$, there is a line bundle $L_x$ on $X$ depending on $x$ such that $T_x^*(E) \cong E \otimes L_x$ for all $x\in X$, where $T_x$ is the translation morphism given by $x\in X$.
\begin{corl}\label{corl4.6}
 Let $E$ be a semi-homogeneous vector bundle of rank $r$ on an abelian variety $X$. Then $\Nef^1\bigl(\Gr(k,E)\bigr) = \overline{\Eff}^1\bigl(\Gr(k,E)\bigr)$ for all integers $k$ with $1\leq k \leq r$.
 \begin{proof}
  By a result due to Mukai, $E$ is Gieseker semistable (see Ch1 \cite{HL10} for definition) with respect to some polarization and it has projective Chern classes zero, i.e., if $c(E)$ is the total Chern class, then $c(E) = \Bigl\{ 1+ c_1(E)/r\Bigr\}^r$ (see p. 260, Theorem 5.8 \cite{Muk78}, p. 266, Proposition 6.13 \cite{Muk78}; also see p. 2 \cite{MN84}). Gieseker semistablity implies slope semistability (see \cite{HL10}). So, in particular, we have $E$ is slope semistable with $c_2\bigl(\End(E)\bigr) = 2rc_2(E) - (r -1)c_1^2(E) = 0$. Hence the result follows.
 \end{proof}
\end{corl}
\begin{corl}\label{corl4.7}
 Let $X$ be a smooth complex projective variety $X$ with $\overline{\Eff}^1(X) = \Nef^1(X)$ and $E_1,E_2,\cdots,E_l$ be finitely many semistable vector bundles on $X$ of ranks $r_1,r_2,\cdots,r_l$ respectively with $c_2\bigl(\End(E_i)\bigr) = 0$ for all $i \in \{1,2,\cdots,l\}$. Fix some integers $k_i$ with $1\leq k_i \leq r_i$ for all $i \in \{
 1,2,\cdots,l\}$. Then we have
  \begin{center}
  $\Nef^1\bigl(\Gr(k_1,E_1) \times_X \Gr(k_2,E_2)\times_X\cdots\times_X \Gr(k_l,E_l)\bigr)$
  
  $= \overline{\Eff}^1\bigl(\Gr(k_1,E_1) \times_X \Gr(k_2,E_2)\times_X\cdots\times_X \Gr(k_l,E_l)\bigr).$
 \end{center}
 
  Moreover, if $\Nef^1(X) = \overline{\Eff}^1(X)$ is a finite polyhedron and
\begin{align*}
\pi_i : X_i := \Gr(k_1,E_1) \times_X \cdots\times_X \Gr(k_i,E_i) \longrightarrow X_{i-1} := \Gr(k_1,E_1) \times_X \cdots\times_X 
\Gr(k_{i-1}, E_{i-1})
\end{align*}
is the projection map for $i\in\{1,2,\cdots,l\}$ with $X_0 = X$, then for each $i$ satisfying $1 \leq i \leq l$
\begin{align*}
\Nef^1\bigl(X_i\bigr) = \overline{\Eff}^1\bigl(X_i\bigr)
= \Bigl\{ x_0\lambda_{\psi_{i}^*E_i,k_i} + x_1\pi_i^*L_1 + x_2\pi_i^*L_2 + \cdots+x_i\pi_i^*L_i\mid x_i\in \mathbb{R}_{\geq 0}\Bigr\}
\end{align*}
where $L_1,L_2,\cdots,L_i$ are the nef generators of the nef cone of $X_{i-1}$ and $\psi_i := \pi_1\circ\pi_2\circ\cdots\circ\pi_{i-1}$.

Conversely, if $X$ is a smooth curve, then the equality 

\begin{center}
$\Nef^1\bigl(\Gr(k_1,E_1) \times_X \cdots\times_X \Gr(k_l,E_l)\bigr)
= \overline{\Eff}^1\bigl(\Gr(k_1,E_1) \times_X \cdots\times_X \Gr(k_l,E_l)\bigr)$ 
\end{center}
implies that $E_i$ is semistable for each $i$.
 \begin{proof}
 We will proceed by induction on $l$. For $l=1$, this is precisely the statement of Theorem \ref{thm4.1}. Now suppose the theorem holds true for $l-1$ many vector bundles. Consider the following fibre product diagram.
 \begin{center}
 \begin{tikzcd}
 X_l:= \Gr(k_1,E_1) \times_X \Gr(k_2,E_2)\times_X \cdots\times_X \Gr(k_l,E_l) \arrow[r, " "] \arrow[d, " \pi_l"]
& \Gr(k_l,E_l)\arrow[d,""]\\
X_{l-1}:= \Gr(k_1,E_1) \times_X \Gr(k_2,E_2)\times_X \cdots\times_X \Gr(k_{l-1}, E_{l-1})  \arrow[r, "\psi_l" ]
& X
\end{tikzcd}
\end{center}
Note that 
\begin{center}
$X_l = \Gr(k_1,E_1) \times_X \Gr(k_2,E_2)\times_X \cdots\times_X \Gr(k_l,E_l) \cong \Gr_{X_{l-1}}\bigl(k_l,\psi_l^*E_l\bigr)$. 
\end{center}
Since $E_l$ is semistable with $c_2\bigl(\End(E_l)\bigr) = 0$ on $X$,  it'        s pullback $\psi_l^*E_l$ under $\psi_l$ is also semistable with $c_2\bigl(\End(\psi_l^*E_l)\bigr) = 0$. By induction hypothesis we have 
\begin{align*}
 \Nef^1\bigl(\Gr(k_1,E_1) \times_X \cdots\times_X \Gr(k_{l-1},E_{l-1})\bigr) = \overline{\Eff}^1\bigl(\Gr(k_1,E_1) \times_X \cdots\times_X \Gr(k_{l-1},E_{l-1})\bigr).
\end{align*}
Therefore applying  Theorem \ref{thm4.1} we get the result. 

Conversely, if $X$ is a curve and $\Nef^1\bigl(\Gr(k_1,E_1) \times_X  \cdots\times_X\Gr(k_l,E_l)\bigr) = \overline{\Eff}^1\bigl(\Gr(k_1,E_1) \times_X \cdots\times_X \Gr(k_l,E_l)\bigr),$  then inductively 
$\Nef^1\bigl(\Gr(k_i,E_i)\bigr) = \overline{\Eff}^1\bigl(\Gr(k_i,E_i)\bigr)$ for each $i\in\{1,2,\cdots,l\}$. This implies that each $E_i$ is semistable by [Corollary 4.3, \cite{B-H-P}]. This completes the proof.
 \end{proof}
 \end{corl}
 \begin{thm}\label{thm4.8}
 Let $E$ and $E'$ be two vector bundles of ranks $m$ and $n$ respectively on a smooth complex projective curve $C$. Fix integers $k$ and $k'$ such that $1\leq k \leq m$ and $1\leq k'\leq n$. Consider the following fibre product diagram:
\begin{center}
 \begin{tikzcd} 
 X = \Gr(k,\pi'^*E) = \Gr(k,E) \times_C \Gr(k',E') \arrow[r, "\pi_1"] \arrow[d, "\pi_2"]
& G = \Gr(k,E)\arrow[d,"\pi"]\\
G' = \Gr(k',E') \arrow[r, "\pi'" ]
& C
\end{tikzcd}
\end{center}
Let $f_1$ and $f_2$ denote the numerical classes of fibres of $\pi$ and $\pi'$ respectively, and $\xi = \pi_1^*c_1\bigl(\mathcal{O}_{\Gr(k,E)}(1)\bigr)$, \hspace{1mm} $\eta = \pi_2^*c_1\bigl(\mathcal{O}_{\Gr(k',E')}(1)\bigr)$, \hspace{1mm} $F = \pi_1^*f_1 = \pi_2^*f_2$.
Then $$\overline{\Eff}^1\bigl(\Gr(k,E) \times_C \Gr(k',E')\bigr) = \Bigl\{a(\xi-\zeta_{E,k}F)+b(\eta-\zeta_{E',k'}F)+cF\mid a,b,c\in \mathbb{R}_{\geq 0}\Bigr\}$$
\begin{proof}
 Note that $\overline{\Eff}^1\bigl(\Gr(k,E)\bigr) = \Bigl\{a(c_1(\mathcal{O}_{\Gr(k,E)}(1))-\zeta_{E,k}f_1)+bf_2\mid a,b\in \mathbb{R}_{\geq 0}\Bigr\}$, and 
 
 $\overline{\Eff}^1\bigl(\Gr(k',E')\bigr) = \Bigl\{a(c_1(\mathcal{O}_{\Gr(k',E')}(1))-\zeta_{E',k'}f_2)+bf_2\mid a,b\in \mathbb{R}_{\geq 0}\Bigr\}$. 
  
  Since the projections $\pi_1$ and $\pi_2$ are smooth maps, the pullback morphisms $\pi_1^*$ and $\pi_2^*$ exist, and hence we have 
  $$\Bigl\{a(\xi-\zeta_{E,k}F)+b(\eta-\zeta_{E',k'}F)+cF\mid a,b,c\in \mathbb{R}_{\geq 0}\Bigr\} \subseteq \overline{\Eff}^1\bigl(\Gr(k,E) \times_C \Gr(k',E')\bigr).$$
  
  It is now enough to show that the generators of the LHS are on the boundary of $\overline{\Eff}^1(X)$.
  
  First note that $$\overline{\Mov}_1(\Gr(k,E)) = \Bigl\{a(\xi_1^{k(r-k)-1}f_1)+b(\xi_1^{k(r-k)}-m_c\xi_1^{k(r-k)-1}f_1\mid a,b\in \mathbb{R}_{\geq0}\Big\},$$
  where $\xi_1 = [\mathcal{O}_{\Gr(k,E)}(1)]$ and $m_c = \deg(\wedge^kE)-\zeta_{E,k}.$
  \vspace{2mm}
  
  We also have ${\pi_1}^*\bigl(\overline{\Mov}_1(G)\bigr) \subset \overline{\Mov}_{\dim(G')}(X)$. Now for any irreducible curve $C'\subset G'$ which is a complete intersection of nef divisors, we have the numerical class $[\pi_2^*(C')] = [G'\times_C C'] \in \Nef^{\dim(G')-1}(X)$. Thus $[\pi_2^*(C')]\cdot \pi_1^*\bigl(\overline{\Mov}_1(G)\bigr) \subset \overline{\Mov}_1(X)$.
  
  Now $$\pi_2^*([C'])\cdot\pi_1^*\bigl(\xi_1^{k(r-k)}-m_c\xi_1^{k(r-k)-1}\bigr)\cdot \pi_1^*(\xi-\zeta_{E,k}F) = 0.$$ This shows that $\xi-\zeta_{E,k}F$ is in the boundary of $\overline{\Eff}^1(X)$. A similar argument will show that $\eta-\zeta_{E',k'}F$ is in the boundary of $\overline{\Eff}^1(X)$. Also $F^2 =0$. 
  
  Therefore,
  $$ \overline{\Eff}^1\bigl(\Gr(k,E) \times_C \Gr(k',E')\bigr) = \Bigl\{a(\xi-\zeta_{E,k}F)+b(\eta-\zeta_{E',k'}F)+cF\mid a,b,c\in \mathbb{R}_{\geq 0}\Bigr\} .
  $$
  \end{proof}
  \end{thm}
  \begin{xrem}
  \rm In the above Theorem \ref{thm4.8}, since $\pi_1$ is a smooth map, in particular it is flat and hence $\pi_1^*$ is an exact functor. Also we observe that for any semistable vector bundle $V$ on $C$ and for any ample line bundle $H$ on $\Gr(k',E')$, ${\pi'}^*(V)$ is $H$-semistable with $\mu_H\bigl({\pi'}^*V) = \mu(V)\bigl(f_2\cdot H^{n-1}\bigr)$ and $f_2\cdot H^{n-1} > 0$. These observations immediately imply that 
\begin{align*}
 0 = {\pi'}^*E_{0} \subset {\pi'}^*E_{1} \subset \cdots\subset {\pi'}^*E_{l-1} \subset {\pi'}^*E_{l} = {\pi'}^*E
\end{align*}
is the unique Harder-Narasimhan filtration of $\pi_1^*E$ with respect to any ample line bundle $H$ on $\Gr(k',E')$ such that each succesive quotient ${\pi'}^*(E_{i}/E_{i-1})$ is semistable with $c_2(\End({\pi'}^*(E_{i}/E_{i-1}))) = 0 \in H^4(\Gr(k',E'),\mathbb{Q})$.
  \end{xrem}

\section{Nef cone of Grassmann Bundle}
We quickly recall our set up. Let $E$ be a vector bundle of rank $ r\geq 2$ on a smooth complex projective variety $X$, and $k$ be an integer such that $1\leq k \leq r$. Let  $\pi_k:\Gr(k,E)\longrightarrow X$ be the projection.
In this section, we give a description of the closed cone $\overline{\NE}\bigl(\Gr(k,E)\bigr)$ of curves in $\Gr(k,E)$  in terms of the closed cone of curves $\overline{\NE}\bigl(\Gr(k,E\vert_C)\bigr)$ for every irreducible curve $C$ in $X$. 
\begin{thm}\label{thm5.1}
 Let $E$ be a  vector bundle of rank $r \geq 2$ on a smooth complex projective variety $X$. For an irreducible curve $C$ in $X$ together with its normalization $\eta_c : \tilde{C}\longrightarrow C$, consider the following fibre product diagram: 
\begin{center}
 \begin{tikzcd}
\tilde{C}\times_X \Gr(k,E)= Gr\bigl(k,\eta_c^*(E\vert_C)\bigr) \arrow[r, "\tilde{\eta_c}"] \arrow[d, "\tilde{\pi_c}"] & C\times_X \Gr(k,E) = \Gr(k,E\vert_C) \arrow[r, "j"] \arrow[d,"\pi_c"] & \Gr(k,E) \arrow[d,"\pi_k"] \\
\tilde{C} \arrow[r, "\eta_c"]          & C \arrow[r, "i"]                                         & X
\end{tikzcd}
\end{center}
where $i$ and $j$ are inclusions. We define $\psi_c:= j \circ \tilde{\eta_c}$. 
Then the closed cone of curves in $\Gr(k,E)$ is given by
 \begin{center}
$ \overline{\NE}\bigl(\Gr(k,E)\bigr) = \overline{\sum\limits_{C\in \Gamma}\bigl(\psi_{c}\bigr)_* \Bigl(\overline{\NE}\bigl(\Gr(k,\eta_c^*(E\vert_C))\bigr)\Bigr)},$
  \end{center}
 where $\Gamma$ is the set of all irreducible curves in $X$.
\end{thm}
\begin{proof}
 Let $\xi \in N^1\bigl(\Gr(k,E)\bigr)_{\mathbb{R}}$ be the numerical equivalence class of the tautological line bundle $\mathcal{O}_{\Gr(k,E)}(1)$ on $\Gr(k,E)$. 
 

 For an irreducible curve $C$ in $X$, we fix the notations $\xi_{\tilde{c}}$ and $f_{\tilde{c}}$ for the numerical equivalence classes of the tautological line bundle 
 $\mathcal{O}_{\Gr(k,\eta_c^*(E\vert_C))}(1)$ of $\Gr\bigl(k,\eta_c^*(E\vert_C)\bigr)$ and a fibre of the map 
 $\tilde{\pi_{c}}$ respectively. 
  The pushforward map  
 $$(\psi_c)_* : N_1\bigl(\Gr(k,\eta_c^*(E\vert_C)\bigr)_{\mathbb{R}} \longrightarrow N_1\bigl(\Gr(k,E)\bigr)_{\mathbb{R}}.$$
 induces the following :
 \begin{align*}
 \overline{\sum\limits_{C\in \Gamma}(\psi_{c})_{*} \Bigl(\overline{\NE}\bigl(\Gr(k,\eta_c^*(E\vert_C))\bigr)\Bigr)} \subseteq \overline{\NE}\bigl(\Gr(k,E)\bigr),
\end{align*}
where the sum is taken over the set $\Gamma$ of all irreducible curves in $X$.
Now to prove the reverse inequality, we consider the  numerical equivalence class $\bigl[\bar{C}\bigr]\in \NE\bigl(\Gr(k,E)\bigr)$ of an irreducible curve $\bar{C}$ in $\Gr(k,E)$ which is not contained in any fibre of $\pi_k$.
Denote $\pi_k(\bar{C})=C$. Then, $\bar{C}\subseteq \Gr(k,E\vert_C)$. 
Then there exists a unique irreducible curve $C'\subseteq \Gr\bigl(k,\eta_c^*(E\vert_C)\bigr)$  such that $\tilde{\eta_c}(C') = \bar{C}$ and $\bigl(\psi_{c}\bigr)_*\bigl([C']\bigr)=\bigl[\bar{C}\bigr]$. Also, the numerical equivalence classes of curves in a fibre of $\tilde{\pi_{c}}$ maps to the numerical classes of curves in a fibre of $\pi_k$ by 
 $\bigl(\psi_{c}\bigr)_*$.
Hence, we have
\begin{align*}
 \overline{\NE}\bigl(\Gr(k,E)\bigr)\subseteq\overline{\sum\limits_{C\in \Gamma}\bigl(\psi_{c}\bigr)_* \Bigl(\overline{\NE}\bigl(\Gr(k,\eta_c^*(E\vert_C))\bigr)\Bigr)}.
\end{align*}
This completes the proof.
\end{proof}
\begin{corl}\label{corl5.2}
Let $X$ be a smooth complex projective surface with
\begin{center} 
$\overline{\NE}(X) = \Bigl\{ a_1[C_1] + a_2[C_2] +\cdots+a_n[C_n] \mid a_i \in \mathbb{R}_{\geq 0}\Bigr\}$
\end{center}
for  some irreducible curves $C_1,C_2,\cdots,C_n$  in $X$. If $E$ is a semistable vector bundle of rank $r\geq2$ on $X$ with $c_2\bigl(\End(E)\bigr) := 2rc_2(E)-(r-1)c_1^2(E) = 0 \in H^4(X,\mathbb{Q})$, then for any integer $k$ with $1\leq k \leq r$, we have
\begin{center}
 $\Nef^1\bigl(\Gr(k,E)\bigr)
   =\Bigl\{ y_0\xi+\pi_k^*\gamma \mid \gamma\in N^1(X)_{\mathbb{R}},  y_0\geq 0,    ky_0\mu(E\vert_{C_j}) + (\gamma\cdot C_j)\geq0 \text{ for all } 1\leq j\leq n\Bigr\}$,
\end{center}
where $\xi$  denotes the numerical equivalence class of 
the tautological bundle $\mathcal{O}_{\Gr(k,E)}(1)$.
\end{corl}
\begin{proof}
Let $C$ be an irreducible curve in $X$ such that $[C] = \sum\limits_{i=1}^n x_{i}[C_{i}] \in N_1(X)_{\mathbb{R}}$ for some $x_i \in \mathbb{R}_{\geq 0}$. As $E$ is semistable with vanishing discriminant, applying Theorem 1.2 in \cite{B-B} to the map
$i \circ \eta_c : \tilde{C} \longrightarrow  C  \hookrightarrow X$, we get that $\eta_c^*\bigl(E\vert_C\bigr)$ is also semistable bundle on $\tilde{C}$ for any irreducible curve $C$ in $X$. Using Proposition 4.1 \cite{B-P}, we have
\begin{align*}
 \Nef^1\bigl(\Gr(k,\eta_c^*(E\vert_C))\bigr) = \Bigl\{ a\bigl(\xi_{\tilde{c}}-k\mu(\eta_c^*(E\vert_C))f_{\tilde{c}}\bigr) + bf_{\tilde{c}} \mid a,b \in \mathbb{R}_{\geq 0}\Bigr\}.
\end{align*} 
We define $l_c:=  \deg\bigl(\wedge^kE\vert_C)-k\mu(E\vert_C)$.
\vspace{2mm}

Applying duality, we then get
\begin{align*}
\overline{\NE}\bigl(\Gr(k,\eta_c^*(E\vert_C))\bigr) = \Bigl\{a(\xi_{\tilde{c}}^{k(r-k)-1}f_{\tilde{c}}) + b(\xi_{\tilde{c}}^{k(r-k)}-l_c\xi_{\tilde{c}}^{k(r-k)-1}f_{\tilde{c}})\mid a,b \in\mathbb{R}_{\geq 0}\Bigr\}.
\end{align*}

Therefore,we have
\vspace{2mm}

  $\overline{\NE}\bigl(\Gr(k,E)\bigr)=\overline{\sum\limits_{C\in \Gamma}(\psi_{c})_* \Bigl(\overline{\NE}\bigl(\Gr(k,\eta_c^*(E\vert_C))\bigr)\Bigr)}$
  \vspace{2mm}
  
  $=
  \overline{\sum\limits_{C\in \Gamma}\Bigl[ \mathbb{R}_{\geq 0}(\xi^{k(r-k)-1}F) + \mathbb{R}_{\geq 0}( \xi^{k(r-k)}\pi_k^*[C] - l_c\xi^{k(r-k)-1}F)\Bigr]}$,
  \vspace{2mm}
  
  where $F$ is the numerical equivalence class of the fibre of the projection map $\pi_k : \Gr(k,E)\longrightarrow X$.
  Note that 


$l_{c_i} = \deg\bigl(\wedge^kE\vert_{C_i})-k\mu(E\vert_{C_i}) = \frac{k}{r}\deg(E\vert_{C_i})\bigl(\rk(\wedge^kE)-1\bigr)$ for every $ 1\leq i \leq n$, and 
$l_c =  \deg\bigl(\wedge^kE\vert_C\bigr)- k\mu(E\vert_C) = \sum\limits_{i=1}^n x_i\{\deg\bigl(\wedge^kE\vert_{C_i}-k\mu(E\vert_{C_i}\bigr)\} = \sum\limits_{i=1}^n x_il_{c_i}.$
\vspace{2mm}

This shows that $\overline{\NE}\bigl(\Gr(k,E)\bigr)$ is generated by
\begin{center}  
  $\Bigl\{ \xi^{k(r-k)-1}F$ , $\xi^{k(r-k)}\pi_k^*[C_i] - l_{c_i}\xi^{k(r-k)-1}F \mid 1\leq i\leq n\Bigr\}$,
\end{center}

and thus applying duality we get the nef cone
\begin{align*}
 \Nef^1\bigl(\Gr(k,E)\bigr)
   =\Bigl\{ y_0\xi+\pi_k^*\gamma \mid \gamma\in N^1(X)_{\mathbb{R}},  y_0\geq 0,    ky_0\mu(E\vert_{C_j}) + (\gamma\cdot C_j)\geq0 \text{ for all } 1\leq j\leq n\Bigr\}.
\end{align*}
\end{proof}
\begin{xrem}\label{xrem2}
  \rm Examples of surfaces $X$ with $\overline{\NE}(X)$ a finite polyhedra include smooth surfaces with Picard number 1, ruled surfaces $\mathbb{P}_C(W)$ over a smooth curve $C$ with $\mu_{\max}(W) = 0$, certain blow-ups $\Bl_x\mathbb{P}_C(W)$ of ruled surfaces $\mathbb{P}_C(W)$ with $\mu_{\max}(W) = 0$, del Pezzo surfaces etc. For a semistable vector bundle $E$ of rank $r$ with $c_2\bigl(\End(E)\bigr)=0$ on such surfaces $X$, one can compute $\Nef^1(\Gr(k,E))$ for every $k$ with $1\leq k \leq r$ using Theorem \ref{thm5.1}. See for Example \ref{exm5.4}.
\end{xrem}
\subsection{The completely decomposable case}\label{direct_sum}
Let $$0 = E_0\subsetneq E_1 \subsetneq E_2 \subsetneq E_3\subsetneq \cdots E_{l-1} \subsetneq E_{l} = E$$ be the Harder-Narasimhan filtration of a vector bundle $E$  of rank $r$ with respect to an ample line bundle $H$ on $X$. We define for any $k$ with $1\leq k\leq r$ the following
$$\theta^H_{E,k} := \bigl(k-\rk(E/E_t)\bigr)\mu_H(E_t/E_{t-1})+\deg_H(E/E_t),$$
where either $t = l$ or $t$ is the smallest integer with $\sum\limits_{i=t+1}^l\rk(E_i/E_{i-1}) = \rk(E/E_t) < k.$

We will use the above notations in what follows.
\begin{corl}\label{cor5.3}
Let $X$ be a smooth complex projective surface with Picard number 1 and $L_X$ be an ample generator of the real N\'{e}ron Severi group $N^1(X)_{\mathbb{R}}$. 
 Let $E = M_1\oplus M_2\oplus\cdots\oplus M_r$ be a completely decomposable vector bundle of rank $r\geq2$ on $X$. 
  We fix an integer $k$ with $1\leq k \leq r$.
  Then
 \begin{center}
   $\Nef^1\bigl(\Gr(k,E)\bigr) = \Bigl\{ y_0\xi + y_1\pi^*L_X \mid y_0\geq 0, y_0\theta^{L_X}_{E,k}+y_1L_X^2 \geq 0\Bigr\}$
  \end{center}
 In particular, if $X = \mathbb{P}^2$, and $E = \mathcal{O}_{\mathbb{P}^2}(a_1)\oplus \mathcal{O}_{\mathbb{P}^2}(a_2) \oplus\cdots\oplus \mathcal{O}_{\mathbb{P}^2}(a_r)$
 a completely decomposable bundle on $\mathbb{P}^2$, then 
 \begin{align*}
 \Nef^1\bigl(\Gr(k,E)\bigr) = \Bigl\{ y_0\xi + y_1\pi^*H \mid y_0\geq 0, y_0\theta_{E,k}^{H} + y_1\geq 0\Bigr\},
 \end{align*} 
 where $H$ is the numerical equivalence class of $\mathcal{O}_{\mathbb{P}^2}(1)$ in $N^1(\mathbb{P}^2)_{\mathbb{R}}$.
 \begin{proof}
 Let $m$ be the least positive integer such that $H^0(mL_X)\neq 0$ and $C_0 \in \vert mL_X \vert$ be an irreducible curve. Then we have $$\overline{\NE}(X) \subseteq \mathbb{R}_{\geq 0}L_X = \mathbb{R}_{\geq 0}\frac{1}{m}[C_0] \subseteq \Nef(X).$$ 
 
 Hence $\overline{\NE}(X) = \mathbb{R}_{\geq 0}[C_0]$.
 Let $C$ be an irreducible curve in $X$, and $[C] \equiv x[C_0]$ for some $x\in \mathbb{R}_{\geq 0}$.
 
 Note that $\theta_{E\vert_C,k} = x\theta_{E\vert_{C_0},k}$ for any $k$ with $1\leq k\leq r$. Also $m\theta^{L_X}_{E,k} = \theta_{E\vert_C,k}$.
 
 Define $l_{C_0} := \deg(\wedge^kE\vert_C) - \theta_{E\vert_{C_0},k}$.
 
 We then have $\overline{\NE}(\Gr(k,E))$ is generated by $\Bigl\{\xi^{k(r-k)-1}F, \xi^{k(r-k)}\pi_k^*[C_0]-l_{C_0}\xi^{k(r-k)-1}F\Bigr\}$.
  
  Therefore applying duality we have 
  \begin{center}
   $\Nef^1\bigl(\Gr(k,E)\bigr) = \Bigl\{y_0\xi + y_1\pi^*L_X \mid y_0 \geq 0, y_0\theta_{E\vert_{C_0},k}+ (y_1L_X\cdot C_0)\geq 0 \Bigr\}$
   
   $= \Bigl\{ y_0\xi + y_1\pi^*L_X \mid y_0\geq 0, y_0\theta^{L_X}_{E,k}+y_1 L_X^2\geq 0\Bigr\}.$
  \end{center}
\end{proof}
 \end{corl}
\begin{exm}\label{exm5.4}
\rm Let $\rho : X = \mathbb{P}(W) \longrightarrow C$ be a ruled surface over a smooth elliptic curve $C$ defined by the rank two bundle 
$W = \mathcal{O}_C \oplus \mathcal{O}_C$. Then $\overline{\NE}(X) = \{ a\zeta + bf \mid a,b\in \mathbb{R}_{\geq 0}\}$, where 
$\zeta = \bigl[\mathcal{O}_{\mathbb{P}(W)}(1)\bigr] \in N^1(X)$ and $f$ is the numerical equivalence class of a fibre of $\rho$.

Let $\rho_x : \tilde{X}_x  =  \Bl_x(X) \longrightarrow X$ be the blow up of $X$ at a closed point $x$ in a section $\sigma$ such that $\mathcal{O}_X(\sigma)\cong \mathcal{O}_{\mathbb{P}(W)}(1)$, and $E_x$ be the exceptional divisor. Then we claim that
$\overline{\NE}(\tilde{X}_x) =\bigl\{ a[C_1] + b[C_2] + c[C_3] \mid a,b,c \in \mathbb{R}_{\geq 0}\bigr\}$, where
$[C_1]=\rho^*_xf - E_x$, $[C_2]=\rho^*_x\zeta - E_x$ and $ [C_3]= E_x$. 

To prove our claim enough to show that if $V$ is an irreducible curve in $X$ which is not a fibre i.e. $[V]=a\zeta+b f$; $a> 0$ and $b\geq 0$, then $\text{mult}_{x}V=r\leq a$.
Let $F$ be a fibre passing through $x$. Then $\text{mult}_{x}F=1$. 

Hence
 $$a=V\cdot F=\sum_{P\in V\cap F}(V.F)_P\geq (V.F)_x\geq r \text{ (see \cite{Har}, Chapter V, Proposition 1.4) }$$

Now we also have the following intersection products:

$b_{11} := C_1\cdot C_1 = (\rho^*_xf - E_x)\cdot(\rho^*_xf - E_x) = E_x^2 = -1$ ; 

$b_{22} := C_2\cdot C_2 =(\rho^*_x\zeta - E_x)\cdot(\rho^*_x\zeta - E_x) = \deg(W) - 1 = -1$ ;
$b_{33} := C_3\cdot C_3 = E_x^2 = -1$.

$b_{12} = b_{21} := C_1\cdot C_2=(\rho^*_xf - E_x)\cdot(\rho^*_x\zeta - E_x) = \zeta\cdot f+ E_x^2 = 0$ ; 

$b_{13} = b_{31} := C_1\cdot C_3 = (\rho^*_xf - E_x)\cdot E_x = 1$ ;
$b_{23} = b_{32} := C_2\cdot C_3 = (\rho^*_x\zeta - E_x)\cdot E_x=1$ ;

Let $E = \rho^*_x\bigl(\rho^*(V)\bigr)$, where $V$ is an indecomposable semistable vector bundle of rank $r$ of degree $d$  on the elliptic curve  $C$. Then $E$ is also semistable bundle of rank $r$ with $c_2\bigl(\End(E)\bigr) = 0$.
Note that, in this example, the Picard number of $\Gr(k,E)$ is 4 for every $k \in\{1,2,\cdots,r\}$.

$\deg(E\vert_{C_1}) = c_1\bigl(\rho^*_x\bigl(\rho^*(V)\bigr)\bigr)\cdot (\rho^*_xf - E_x) = 0$ ;

$\deg(E\vert_{C_2}) = c_1\bigl(\rho^*_x\bigl(\rho^*(V)\bigr)\bigr)\cdot (\rho^*_x\zeta - E_x) = \deg(V) = d$;

$\deg(E\vert_{C_3}) = c_1\bigl(\rho^*_x\bigl(\rho^*(V)\bigr)\bigr)\cdot  E_x = 0$.

$y_0\mu(E\vert_{C_1}) + y_1b_{11} + y_2b_{21} +y_3b_{31} = y_3 -y_1$ ;

$y_0\mu(E\vert_{C_2}) + y_1b_{12} + y_2b_{22} +y_3b_{32} = y_0\mu(V) + y_3 -y_2$ ;

$y_0\mu(E\vert_{C_3}) + y_1b_{13} + y_2b_{23} +y_3b_{33} = y_1 +y_2 -y_3$.
\vspace{1mm}

Therefore, the nef cone of the grassmann bundle $\pi : \Gr(k,E) \longrightarrow \tilde{X}_x$ for $k\in \{1,2,\cdots,r\}$ is 
\begin{align*}
 \Nef^1\bigl(\Gr(k,E)\bigr) = \Bigl\{ y_0\xi+\sum\limits_{i=1}^3y_i\pi^*(C_i) \mid  y_0\geq 0, y_3 -y_1\geq 0, ky_0\mu(V) + y_3 -y_2\geq 0,  y_1 +y_2 -y_3\geq 0\Bigr\}.
\end{align*}
\end{exm}

\begin{thm}\label{thm5.5}
 Let $E$ and $E'$ be two vector bundles of ranks $m$ and $n$ respectively on a smooth complex projective curve $C$. Fix integers $k$ and $k'$ such that $1\leq k \leq m$ and $1\leq k'\leq n$. Consider the following fibre product diagram:
\begin{center}
 \begin{tikzcd} 
 X = \Gr(k,\pi'^*E) = \Gr(k,E) \times_C \Gr(k',E') \arrow[r, "\pi_1"] \arrow[d, "\pi_2"]
& G = \Gr(k,E)\arrow[d,"\pi"]\\
G' = \Gr(k',E') \arrow[r, "\pi'" ]
& C
\end{tikzcd}
\end{center}
Let $f_1$ and $f_2$ denote the numerical classes of fibres of $\pi$ and $\pi'$ respectively, and $\xi = \pi_1^*c_1\bigl(\mathcal{O}_{\Gr(k,E)}(1)\bigr)$, \hspace{1mm} $\eta = \pi_2^*c_1\bigl(\mathcal{O}_{\Gr(k',E')}(1)\bigr)$, \hspace{1mm} $F = \pi_1^*f_1 = \pi_2^*f_2$.
Then $$\Nef^1\bigl(\Gr(k,E) \times_C \Gr(k',E')\bigr) = \Bigl\{a(\xi-\theta_{E,k}F)+b(\eta-\theta_{E',k'}F)+cF\mid a,b,c\in \mathbb{R}_{\geq 0}\Bigr\}$$
\begin{proof}
 Note that $\Nef^1\bigl(\Gr(k,E)\bigr) = \Bigl\{a(c_1(\mathcal{O}_{\Gr(k,E)}(1))-\theta_{E,k}f_1)+bf_2\mid a,b\in \mathbb{R}_{\geq 0}\Bigr\}$, and 
 
 $\Nef^1\bigl(\Gr(k',E')\bigr) = \Bigl\{a(c_1(\mathcal{O}_{\Gr(k',E')}(1))-\theta_{E',k'}f_2)+bf_2\mid a,b\in \mathbb{R}_{\geq 0}\Bigr\}$. 
  
  Since pullback of nef divisors are nef, we have  
  $$\Bigl\{a(\xi-\theta_{E,k}F)+b(\eta-\theta_{E',k'}F)+cF\mid a,b,c\in \mathbb{R}_{\geq 0}\Bigr\} \subseteq \Nef^1\bigl(\Gr(k,E) \times_C \Gr(k',E')\bigr).$$
  
  To check the reverse inclusion, enough to show that the generators of the LHS form the boundary of $\Nef^1\bigl(\Gr(k,E) \times_C \Gr(k',E')\bigr)$.  
  
  Let $\gamma\in \overline{\NE}(G)\subset N_1(G)$. Then $\pi_1^*(\gamma)\in N_{\dim(G)}(X)\cong N^{\dim(X)-\dim(G')}(X)\cong N^{\dim(G)-1}(X)$. If $C'\subset G'$ be an irreducible curve in $G'$ which is a complete intersection of nef divisors, then note that $[G\times_C C']\in N_{\dim(G)}(X)$. Now clearly $\pi_1^*\gamma\cdot [G\times_C C']\in N_1(X)$  is numerical class of an effective 1-cycle in $X$.  Observe that $\pi_1^*\gamma\cdot [G\times_C C']$ is a class of effective curve in $X$ and whose image by the projection $\pi_1$ is an effective curve isomorphic to an element of curve class $\gamma$. 
  
  We recall that
  $$\overline{\NE}(G)=\Bigl\{a\xi_1^{k(r-k)-1}f_1+b(\xi_1^{k(r-k)}-l_c\xi_1^{k(r-k)-1}f_1)\mid a,b\in\mathbb{R}_{\geq0}\Bigr\},$$
where $\xi_1 \equiv [\mathcal{O}_{\Gr(k,E)}(1)]$ and $l_c=\deg(\wedge^k E)- \theta_{E,k}$.
If we take $\gamma= \xi_1^{k(r-k)}-l_c\xi^{k(r-k)-1}f_1$, then 
$$(\xi-\theta_{E,k}F)\cdot \pi_1^*\gamma\cdot [G\times_C C']=\pi_1^*((\xi_1-\theta_{E,k}f_1)\cdot \gamma)\cdot [G\times_C C']=0.$$

Hence $\xi-\theta_{E,k}F$ is in the boundary of $\Nef^1(X)$. Similar arguments hold for other generator, and $F^2 = 0$. Hence the result follows.
\end{proof}
\end{thm}
\section{Ampleness of universal quotient bundles}
\begin{thm}\label{thm6.1}
  Let $E$ be a vector bundle of rank $r$ on a smooth complex projective variety $X$. For every $k$ with $1\leq k\leq r$, consider the universal quotient bundle $Q_k$ on the grassmann bundle $\Gr(k,E)$ over $X$. Consider the following commutative diagram:
  \begin{center}
   \begin{tikzcd}
\mathbb{P}_Y(Q_k) \arrow[rd, "\phi_k"] \arrow[r, "\psi_k"] & Y=\Gr(k,E) \arrow[d,"\pi_k"]\\
& X
\end{tikzcd}
  \end{center}
Then the following are equivalent:
  \begin{enumerate}
   \item $\Phi_{E,k} := c_1(\mathcal{O}_{\mathbb{P}(Q_k)}(1)) - \frac{1}{r}\phi_k^*c_1(E)$ is nef for every $k$ with $1\leq k\leq r.$
   \item $E$ is semistable with $c_2\bigl(\End(E)\bigr)=0\in H^4(X,\mathbb{Q})$.
  \end{enumerate}
\begin{proof}
(1)$\implies$(2) Let $f: C\longrightarrow X$ be a nonconstant morphism from a smooth projective curve $C$. Note that for any $k$ with $1\leq k\leq r$, we have $$\Phi_{f^*(E),k} = \tilde{f}^*(\Phi_{E,k}),$$
where $\tilde{f} : \mathbb{P}(Q_k')\longrightarrow \mathbb{P}(Q_k)$ is the map induced by $f$, and $Q_k'$ is the universal quotient bundle of $f^*(E)$.  

As each $\Phi_{E,k}$ is nef, we have each $\Phi_{f^*(E),k}$ is nef. Therefore, we conclude that $f^*(E)$ is semistable on $C$ by [Theorem 1.1,\cite{BH}]. This shows that $E$ is a semistable bundle on $X$ with $c_2(\End(E))=0$.

(2)$\implies$(1) Suppose $\Phi_{E,k}$ is not nef for some $k$. Then there exists an irreducible curve $C'\subset \mathbb{P}(Q_k)$ such that $C'\cdot \Phi_{E,k} <0$. 

Let $h:C''\longrightarrow X$ be a finite morphism from a smooth curve $C''$, and consider the commutative diagram:
\begin{center}
 \begin{tikzcd} 
 \mathbb{P}(Q_k') \arrow[r, "h'"] \arrow[d, "\psi_k'"]
& \mathbb{P}(Q_k) \arrow[d,"\psi_k"]\\
\Gr(k,h^*E) \arrow[r, "\tilde{h}" ] \arrow[d, "\pi_k'"]
& \Gr(k,E) \arrow[d, "\pi_k"] \\
C''\arrow[r, "h"] & X
\end{tikzcd}
\end{center}
We may choose the pair $(C'',h)$ in such a way that the fiber product $\tilde{C} = C''\times_C C'$ (a curve in $\mathbb{P}(Q_k')$) is a union of curves $C_j$ which project onto $C''$ with degree 1 and meet the fiber $f_k'$ of $\pi_k'$ at just one point. 

One has $[C_j]\cdot \Phi_{h^*E,k} <0$. Consider the surjection $\Gamma_j : C_j \longrightarrow \Gr(h^*E,k)$. We denote by $Q_j$ the restriction of $Q_k'$ to $C_j$ under $\Gamma_j$. Now consider the following vector bundle $E_j := \bigl(\tilde{h}\circ \pi_k^*E\bigr)\vert_{\Gamma_j}$. We have an epimorphism $E_j\longrightarrow Q_j\longrightarrow 0$. The composition $h\circ \pi_k'\vert_{\Gamma_j} = h_j$ is a finite morphism so that $E_j$ is semistable. 
Now we get $$\mu(Q_j) \leq \frac{[C_j]\cdot \xi_j}{{\pi_k'}_*([C_j])} = [C_j]\cdot \xi_k' = [C_j]\cdot \bigl(\Phi_{h^*E,k}+\mu(h^*E)f_k'\bigr)<\mu(E_j),$$ which is a contradiction, as $E_j$ is semistable. 
This concludes the proof.
\end{proof}
\end{thm}
\begin{corl}
 Let $E$ be a semistable  nef vector bundle of rank $r$ on a smooth complex projective variety $X$ of dimension $n$ with $c_2\bigl(\End(E)\bigr) = 0 \in H^4(X,\mathbb{Q})$. Then for every $k$ with $1\leq k \leq r$, the universal quotient bundle $Q_k$ is  nef on the grassmann bundle $\Gr(k,E)$.
 \begin{proof}
  Note that $\det(Q_k) = \mathcal{O}_{\Gr(k,E)}(1) = \omega^*\mathcal{O}_{\mathbb{P}(\wedge^kE)}(1)$ where $\omega : \Gr(k,E)\hookrightarrow \mathbb{P}(\wedge^kE)$ is the Plucker embedding. Since $E$ is  nef, we have $\wedge^k(E)$ is also  nef. Thus $\det(Q_k)$ is  nef.  
  Also, by the previous  Theorem \ref{thm6.1},  we have $c_1(\mathcal{O}_{\mathbb{P}(Q_k)}(1)) - \frac{1}{r}\pi_k^*c_1(E)$ is nef. Thus we have $c_1(\mathcal{O}_{\mathbb{P}(Q_k)}(1))$ is nef, i.e. $Q_k$ is nef.
  \end{proof}
\end{corl}
\begin{prop}
 Let $E$ be an ample vector bundle of rank $r$ on a smooth complex projective variety $X$ of dimension $n$.Then for every $k$ with $1\leq k \leq r$, the determinant bundle of the universal quotient bundle $\det(Q_k)$ is ample.
 \begin{proof}
  Note that $\det(Q_k) = \mathcal{O}_{\Gr(k,E)}(1) = \omega^*\mathcal{O}_{\mathbb{P}(\wedge^kE)}(1)$ where $\omega : \Gr(k,E)\hookrightarrow \mathbb{P}(\wedge^kE)$ is the Plucker embedding. Since $E$ is ample, we have $\wedge^k(E)$ is also  ample. As $\omega$ is an embedding, we  thus conclude that $\det(Q_k)$ is ample.
 \end{proof}
\end{prop}
\begin{xrem}
 \rm Note that a semistable bundle $E$ on a smooth curve $C$ is ample if and only if $\det(E)$ is ample. However, analogous statement is not true for semistable bundle on higher dimensional varieties.  Hence $\det(Q_k)$ being ample does not ensure the ampleness of $Q_k$ whenever dimension of $X$ is at least 2.
\end{xrem}
\begin{corl}
 Let $E$ be a semistable vector bundle of rank $r$ on a smooth complex projective surface $X$ of dimension $n$ with $c_2\bigl(\End(E)\bigr) = 0 \in H^4(X,\mathbb{Q})$. Let $H$ be an ample line bundle on $X$. Assume that for some $k$ with $1\leq k \leq r$, we have $\mu_H(Q_k) = \mu_H(\pi^*E)$. Then the universal quotient bundle $Q_k$ is a semistable bundle on the grassmann bundle $\Gr(k,E)$ with $c_2\bigl(\End(Q_k)\bigr)=0 \in H^4(\Gr(k,E),\mathbb{Q})$. Consequently, $$\Nef^1\bigl(\Gr(s,Q_k)\bigr) = \overline{\Eff}^1\bigl(\Gr(s,Q_k)\bigr)$$ for every $s$ such that $1\leq s \leq k$.
 Moreover, in addition, if $E$ is ample, then $Q_k$ is also ample.
 \begin{proof}
  We consider the following exact sequence :
  $$0\longrightarrow S_{r-k}\longrightarrow \pi_k^*E\longrightarrow Q_k\longrightarrow 0.$$
   From the above exact sequence we have 
  $$\frac{r}{r-k}\Delta(S_k) +\frac{r}{k}\Delta(Q_k)  - \frac{\xi^2}{k(r-k)} = \Delta(\pi_k^*E) = 0,$$
  where $\Delta(S_{r-k}) = c_2(\End(S_k)), \Delta(Q_k) = c_2(\End(Q_k))$ and $\xi \equiv  rc_1(S_{r-k})-(r-k)c_1(\pi^*E) $.
  
  As $\mu_H(Q_k) = \mu_H(\pi_k^*E) = \mu_H(S_{r-k})$, and $\pi_k^*E$ is semistable, we  conclude that both $Q_k$ and $S_{r-k}$ are semistable bundles. By Bogomolov inequality, we have $\Delta(S_k)  \geq 0, \Delta(Q_k)\geq 0$. 
  
  Also $\xi \cdot H =0$. So by Hodge index Theorem $\xi^2\leq 0$. Hence we conclude $Q_k$ is semistable with $c_2(\End(Q_k)) = 0$, and consequently  using Theorem \ref{thm4.1} we get $\Nef^1\bigl(\Gr(s,Q_k)\bigr) = \overline{\Eff}^1\bigl(\Gr(s,Q_k)\bigr)$ with $1\leq s \leq k$.
  
  Now, in addition, if $E$ is ample, then $\det(Q_k)$ is ample. Therefore, $Q_k$ is also ample by [Theorem 1, \cite{MR21}] as $Q_k$ is a semistable bundle with $c_2(\End(Q_k)) =0$. This completes the proof.
 \end{proof}
\end{corl}

\section{Acknowledgement}

The authors would like to thank A J Parameswaran and Paramesh Sankaran  for many fruitful discussions. The first author is supported financially by SERB-NPDF fellowship (File no : PDF/2021/00028).


\begin{thebibliography}{************}

\normalsize
\baselineskip=17pt

\bibitem[BB08]{B-B}{Indranil Biswas and Ugo Bruzzo,}
\emph{On Semistable Principal Bundles over a Complex Projective Manifold,}
International Mathematics Research Notices, Article ID rnn035, (2008).

\bibitem[BDPP13]{BDPP13} S. Boucksom, J.P. Demailly, M. P\u{a}un and T Peternell,
\emph{The pseudo-effective cone of a compact K\"{a}hler manifold and varieties of negative Kodaira dimension,}
Journal of Algebraic Geometry. 22(2), (2013) 201-248.

\bibitem[BH06]{BH} U.Bruzzo and D.Hern\'{a}dez Ruip\'{e}rez, 
\emph{Semistability vs.nefness for (Higgs)vector bundles.}
Differential Geometry and its Applications 24 (2006) 403-416.

\bibitem[BHP14]{B-H-P} Indranil Biswas, Amit Hogadi and A.J.Parameswaran,
\emph{Pseudo-effective cone of Grassmann Bundles over a curve,}
Geom Dedicata, 172, (2014) 69-77.

\bibitem[BP14]{B-P} Indranil Biswas, A.J. Parameswaran
\emph{Nef cone of flag bundles over a curve}
Kyoto Journal of Mathematics, Vol 54, No. 2(2014) 353-366.



\bibitem[F11]{Fu} Mihai Fulger,
\emph{The cones of effective cycles on projective bundles over curves,}
Math.Z.,  269, (2011) 449-459.


\bibitem[F98]{F98} William Fulton.
\emph{Intersection Theory,}
Second Edition, Springer, (1998).


\bibitem[H77]{Har} Robin Hartshorne,
\emph{Algebraic Geometry,} 
Graduate Text in Mathematics, Springer, (1977).

\bibitem[HL10]{HL10} Daniel Huybrechts and Manfred Lehn.
\emph{The Geometry of Moduli Spaces of Sheaves,}
Second Edition, Cambridge University Press, (2010).

\bibitem[L07]{L1} Robert Lazarsfeld,
\emph{Positivity in Algebraic Geometry, Volume I,}
Springer, (2007).

\bibitem[L11]{L2} R. K. Lazarsfeld,
\emph{Positivity in Algebraic Geometry I,} 
A Series of Modern Surveys in Mathematics  48  (Springer-Verlag Berlin Heidelberg).


\bibitem[M87]{Miy87} Yoichi Miyaoka, 
\emph{The Chern classes and Kodaira dimension of a minimal variety,}
Algebraic geometry, Sendai,
1985, Adv. Stud. Pure Math., vol. 10, North-Holland, Amsterdam, 1987, pp. 449-476. 

\bibitem[N99]{N99} Noboru Nakayama,
\emph{Normalized Tautological divisors of semi-stable vector bundles,}
(Japanese) Free resolutions of coordinate rings of projective varieties and related topics (Japanese) (Kyoto, 1998). Sūrikaisekikenkyūsho Kōkyūroku No. 1078 (1999), 167-173.

\bibitem[M78]{Muk78} Shigeru Mukai,
\emph{Semi-homogeneous vector bundles on an abelian variety,}
J. Math. Kyoto Univ. (JMKYAZ) 18-2, (1978) 239-272.

\bibitem[MN84]{MN84} V B Mehta and Madhav Nori,
\emph{Semistable sheaves on homogeneous spaces and abelian varieties,}
Proc. Indian Acad. Sci. (Math Sci) Vol 93, No 1, November (1984) pp 1-12.

\bibitem[MR21]{MR21} S. Misra and N. Ray, 
\emph{On Ampleness of vector bundles,}
C. R. Math. Acad. Sci. Paris 359 (2021), 763-772.



\end{thebibliography}
\end{document}